\documentclass[onefignum,onetabnum]{siamonline171218}

\usepackage{mathtools}
\usepackage[]{overpic}
\usepackage{bm}
\usepackage{subcaption}
\usepackage[usestackEOL]{stackengine}
\usepackage{booktabs}
\usepackage{amsmath}

\DeclareMathOperator*{\argmin}{argmin}
\newcommand{\norm}[1]{\left\lVert#1\right\rVert}
\DeclarePairedDelimiter\abs{\lvert}{\rvert}

\usepackage{pifont}
\newcommand{\xmark}{\ding{53}}%

\usepackage{lipsum}
\usepackage{amsfonts}
\usepackage{graphicx}
\usepackage{epstopdf}
\usepackage{algorithmic}
\ifpdf
  \DeclareGraphicsExtensions{.eps,.pdf,.png,.jpg}
\else
  \DeclareGraphicsExtensions{.eps}
\fi

\usepackage{enumitem}
\setlist[enumerate]{leftmargin=.5in}
\setlist[itemize]{leftmargin=.5in}

\newsiamremark{remark}{Remark}
\newsiamremark{hypothesis}{Hypothesis}
\crefname{hypothesis}{Hypothesis}{Hypotheses}
\newsiamthm{claim}{Claim}

\headers{Spatiotemporal Decomposition of Smooth Signals}{S. M. Hirsh, B. W. Brunton, and J. N. Kutz}

\title{Data-driven Spatiotemporal Modal Decomposition for Time Frequency Analysis}

\author{Seth M. Hirsh\thanks{Department of Physics, University of Washington, Seattle, WA 
  (\email{hirshs@uw.edu}).}
\and Bingni W. Brunton\thanks{Department of Biology, University of Washington, Seattle, WA 
  (\email{bbrunton@uw.edu}).}
\and J. Nathan Kutz\footnotemark[3] \thanks{Department of Applied Mathematics, University of Washington, Seattle, WA (\email{kutz@uw.edu}}).}

\usepackage{amsopn}

\ifpdf
\hypersetup{
  pdftitle={Data-driven Spatiotemporal Modal Decomposition for Time Frequency Analysis}
  pdfauthor={S. M. Hirsh, B. W. Brunton, and J. N. Kutz}
}
\fi

\newcommand{\bunderbrace}[2]{%
  \begin{array}[t]{@{}c@{}}
  \underbrace{#1}\\
  #2{}
  \end{array}
}

\DeclarePairedDelimiterX\set[1]\lbrace\rbrace{#1}

\newcommand{\R}{\mathbb{R}}

\begin{document}

\maketitle

\begin{abstract}
We propose a new solution to the blind source separation problem that factors mixed time-series signals into a sum of spatiotemporal modes, with the constraint that the temporal components are intrinsic mode functions (IMF's).
The key motivation is that IMF's allow the computation of meaningful Hilbert transforms of non-stationary data, from which instantaneous time-frequency representations may be derived.
Our spatiotemporal intrinsic mode decomposition (STIMD) method leverages spatial correlations to generalize the extraction of IMF's from one-dimensional signals, commonly performed using the empirical mode decomposition (EMD), to multi-dimensional signals.
Further, this data-driven method enables future-state prediction.
We demonstrate STIMD on several synthetic examples, comparing it to common matrix factorization techniques, namely singular value decomposition (SVD), independent component analysis (ICA), and dynamic mode decomposition (DMD).
We show that STIMD outperforms these methods at reconstruction and extracting interpretable modes.
Next, we apply STIMD to analyze two real-world datasets, gravitational wave data and neural recordings from the rodent hippocampus.
\end{abstract}




\section{Introduction}

The analysis of spatiotemporal signals is of critical importance for characterizing emerging large-scale measurements in wide variety of scientific and engineering applications.
Significant advances in sensor cost, data storage, and processing power has led to a rapidly increasing availability of data in domains including neuroscience, atmospheric sciences, and finance, to name a few.
For the majority of these applications, the underlying dynamics are nonlinear, non-stationary, and the governing equations are poorly known at best.
Therefore, data-driven modeling tools have become increasingly central, and the ability to extract interpretable structure and provide physical insights are crucial for advancing the field.
Towards this goal, here we introduce the {\em spatiotemporal intrinsic mode decomposition} (STIMD) method, which factors spatiotemporal data into a product of spatial modes and temporal modes, with the constraint that the temporal modes are \emph{intrinsic mode functions} (IMF's~\cite{huang1998empirical,Huang_2008}).  
Our method allows us to perform instantaneous time frequency analysis by computing a Hilbert transform of the data; in addition, it is possible to make future-state predictions of the spatiotemporal system.

The Fourier transform is a very widely used technique for analyzing the power spectral features in time-series signals.
However, this technique assumes periodic systems and performs poorly when the signal is nonlinear and non-stationary.
A large variety of very successful methods have been developed to work with such signals, including windows versions of Fourier analysis and wavelet analysis~\cite{daubechies1992ten}.
One related method is the empirical mode decomposition (EMD~\cite{huang1998empirical,rilling2003empirical,wu2007trend,huang2014hilbert}), which had been developed with the motivation to compute instantaneous time frequency analysis of non-stationary signals.
EMD decomposes real-valued signals into a set of intrinsic mode functions (IMF's), which have the feature that they are suitable for computing meaningful Hilbert transforms.
EMD has been widely applied in a variety of application domains~\cite{huang2003applications,pigorini2011time,wang2012comparing}.
Thus, EMD is able to, analyze non-stationary time-series data, where frequencies vary in time.,
However, EMD an empirical algorithm that had been developed without a rigorous mathemati,cal foundation \cite{rilling2003empirical}. 
To provide this foundational theory, several methods such as the synchrosqueezed wavelet transform~\cite{daubechies2011synchrosqueezed} and the nonlinear matching pursuit method (NMP,~\cite{hou2013data,hou2014convergence}) have been developed to provide a rigorous theoretical foundation for the EMD approach. 
 
For systems with spatiotemporal dynamics, it is often possible and desirable to use several sensors placed at different locations to simultaneously gather data about the system.
These additional measurements and their correlations in space may be leveraged to produce more accurate models.
A large family of methods have been developed for factoring spatiotemporal data into products of two sets modes---one spatial and one temporal.
This factorization problem is also known as the blind source separation problem; in other words, the goal is to extract and disambiguate the underlying signals that comprise the measurement data.
The solution to the decomposition is generally underdetermined, and various results can be obtained by making different assumptions~\cite{ghahramani2004unsupervised}.
We give a brief overview of some of the most widely used decomposition methods in Section~\ref{sec:relatedwork}.

In this paper, we propose a new solution to the blind source separation problem for spatiotemporal non-stationary signals.
This \emph{spatiotemporal intrinsic mode decomposition} (STIMD) is motivated by EMD and builds on the NMP method to factor spatiotemporal data into a set of spatial modes and IMF temporal modes.
To our knowledge, no other decomposition has been described to satisfy these assumptions.
With STIMD, we can compute an instantaneous time-frequency representation with a Hilbert transform and also perform future state prediction.
In Section~\ref{sec:method}, we describe the STIMD method and characterize its behavior on several synthetic non-stationary time-series data examples; in particular, we focus on signals containing frequency modulation.
We show that STIMD extracts the underlying source signals more accurately and reliably than several other commonly used factor analysis techniques.
Further, we illustrate its dependence on noise magnitude and initial conditions.
Next, in Section~\ref{sec:demos}, we apply our method to two real-life datasets, namely measurements of gravitational waves from the laser interferometer gravitational observatory (LIGO) experiment and recordings from neural activity from the rodent hippocampus.  
Our results show that leveraging the architectures jointly greatly improves the performance of spatiotemporal decompositions.
A summary and future improvements are found in Section~\ref{sec:conclusions}.  

\section{Related Work}
\label{sec:relatedwork}

Our STIMD algorithm is based on a number of recent innovations in time frequency analysis of single signals, specifically the empirical mode decomposition (EMD) and nonlinear matching pursuit (NMP).
We first describe these techniques and summarize their approach as algorithms in Sections~\ref{subsec:EMD}, \ref{subsec:nmp}, \ref{subsec:nmpperiodic}.

Next, we give an overview of blind source separation and factor analysis methods to decompose spatiotemporal signals in Section~\ref{subsec:factoring}.
We highlight the common structure and differing assumptions of three widely used techniques: singular value decomposition (SVD, Section~\ref{subsec:svd}), independent component analysis (ICA, Section~\ref{subsec:ica}), and dynamic mode decomposition (DMD, Section~\ref{subsec:dmd}).

\subsection{EMD and the Hilbert Transform}
\label{subsec:EMD}
Consider a signal $f(t) : [t_0,t_1] \to \mathbb{R}$ on which we would like to perform time-frequency analysis, extracting both the temporal and frequency features of the signal simultaneously. 
If the signal is stationary, we may choose to perform the Fourier transform, which decomposes the signal into a basis of sines and cosines; this basis may also be called a dictionary \cite{hou2011adaptive}. 

For non-stationary signals, a dictionary of sines and cosines do not well represent the signal, so we must choose to decompose our signals using a different basis set. 
One possible dictionary is the set of all intrinsic mode functions (IMFs); an IMF has the important property that it has a well-defined Hilbert spectrum. 

IMFs are defined by the following criteria:
 
\begin{enumerate}
    \item The number of extrema and the number of zero crossings of the function must be equal (or differ by at most one).
    
    \item At any point of the function, the average of the upper envelope and the lower envelope defined by the local extrema must be zero; in other words, the function is symmetric with respect to zero).
\end{enumerate}

Mathematically, all real-valued functions $s(t)$ obeying these criteria may be expressed in the form
\begin{equation}
	s(t) = a(t) \cos(\theta(t))
\end{equation}
for some $a,\theta : [t_0, t_1] \to \R$ such that $\theta'(t) > 0$.

One notable property of an IMF is that it has a well-defined Hilbert spectrum.  
Specifically, the function $s(t)$ has the analytic continuation $\tilde{s}(t) : \mathbb{R} \to \mathbb{C}:$
\begin{equation}
	\tilde{s}(t) = a(t) e^{i \theta(t)},
\end{equation}
which has a well-defined instantaneous frequency
\begin{equation}
	\omega(t) = \frac{d \theta}{dt}.
\end{equation}
In 1998, Huang introduced a method for decomposing a signal $x(t)$ into a sum of IMFs $s_j(t)$ and a residual $\rho(t)$ through a recursive sifting process known as EMD~\cite{huang1998empirical}:  
\begin{equation}
	x(t) = \sum_j s_j(t) + \rho(t) .
\end{equation}
Briefly, at each recursive step of EMD, a cubic spline is fit to the local minima and maxima of the data, forming two envelopes. 
The mean of these envelopes $m(t)$ is then subtracted from $x(t)$ to form a residual. 
If the residual $\rho(t)$ is an IMF, it is extracted and the process is applied to the remainder of the data. 
Otherwise, this process is applied recursively to the residual until an IMF is obtained.
Specifically, the $j$th IMF computed after $k$ iterations (assuming that it satisfies the definition of an IMF) is 
\begin{equation}
  s_j(t) = x_j(t) - \sum_{i = 1}^k m_{j,i}(t),
\end{equation}
where $x_j(t) = \sum_{l = 1}^{j - 1} s_l(t)$ and $m_{j,i}(t)$ is the mean of the envelopes computed after $j$ iterations \cite{huang1998empirical,bellini2014final}.

This method has been demonstrated to be successful in practice on a wide number of applications~\cite{wang2012comparing,pigorini2011time,huang2003applications}. In addition, several multivariate and multidimensional extenstions have been developed~\cite{rilling2007bivariate,mandic2013empirical,feng2014fast}\footnote{To the best of our knowledge none of these previous methods perform a matrix factorization comparable to the STIMD method described in Section \ref{sec:method}.}. 
However, EMD is empirical in nature and its mathematical foundation is still poorly understood; therefore, the next section describes one mathematical architecture to establish EMD on a rigorous foundation.

\subsection{Nonlinear Matching Pursuit Method (NMP)} \label{subsec:nmp}
One recent innovation for making EMD rigorous is the nonlinear matching pursuit (NMP) method developed by Hou \& Shi~ \cite{hou2013data}.  
Like EMD, NMP decomposes a signal $x(t)$ into a sum of IMF's $s_k(t)$. 
In particular, for NMP we assume that $x(t)$ can be represented by only a few IMF's. 
Thus, the goal of NMP is to solve the optimization problem
\begin{equation}
	\min M  \text{ such that } x(t) = \sum_{j = 1}^M s_j(t).
\end{equation}
We further assume in this algorithm that the IMF's $s_j$ contain only interwave frequency modulation, as defined in the following definition. 
Since $\theta(t)$ is monotonic by the Invertible Function theorem, we can express an IMF $s(\cdot)$ as a function of $\theta$
\begin{equation}
	s(\theta) = a(\theta) \cos(\theta).
\end{equation}
$s(\theta)$ is defined to have interwave frequency modulation if $a(\theta)$ and $\theta'(\theta) := \frac{d \theta}{dt}|_\theta$ are smoother than $\cos(\theta)$\footnote{Since $\theta' > 0$, by the Inverse Function Theorem, there is a one-to-one mapping between $t$ and $\theta$, so that we can thus express the IMF in $\theta$ space without losing information.}. 
Saying that $a(\theta)$ and $\theta'$ are smoother than $\cos(\theta)$ means that $a(\theta)$ are in the set
\begin{equation}
  V(\theta, \lambda) = \text{span} \left\{ 1, \cos\left(\frac{k \theta}{2 L_{\theta}}\right), \sin \left(\frac{k \theta}{2 L_{\theta}} \right) : 1 \leq k \leq 2 \lambda L_{\theta} \right\},
\label{eq:Vtheta}
\end{equation}
where $L_{\theta} = \left\lfloor \frac{\theta(t_1) - \theta(t_0)}{2 \pi} \right\rfloor$ and $\lambda = 1/2$. 
Note that the parameter $\lambda$ is important in the implementation of NMP, as described in Section~\ref{subsec:nmpperiodic}.

Physically, the signals with interwave frequency modulation roughly correspond to solutions of second order differential equations of the form 
\begin{equation}
  \ddot{x} + b(t) \dot{x} + c(t) x = 0,
\end{equation}
where $b(t)$ and $c(t)$ are sufficiently smooth. 
More details about this solution can be found in~\cite{hou2015sparse}.

To summarize, the dictionary of interwave frequency modulated IMF's is 
\begin{equation}
	\mathcal{D} = \set{a(\theta) \cos(\theta) : a \in V(\theta), \theta' \in V(\theta) \text{ and }\theta'(t) \geq 0},
\end{equation}
and the minimization problem for the NMP method becomes
\begin{equation}
    \min_{a_k, \theta_k} M  \text{ such that }  x = \sum_{k = 1}^M a_k \cos{\theta_k} \text{ and } a_k \cos{\theta_k} \in \mathcal{D} \text{   } \forall k.
    \label{eq:NMP}
\end{equation}
In the case of signals with noise, the equality in \eqref{eq:NMP} is replaced with the inequality $\abs{x - \sum_{k = 1}^M a_k \cos{\theta_k} }_2 \leq \delta$.

\subsection{NMP Implementation for Periodic Signals} \label{subsec:nmpperiodic}
The NMP minimization problem is solved using matching pursuit.  
As with the EMD algorithm, each IMF is discovered by a greedy optimization and subtracted from the residual $r_k(t)$ at each step. 
Specifically, given the signal $x(t)$ we extract the first IMF by solving the minimization problem 
\begin{equation}
  \argmin_{a, \theta} \abs{x(t) - a(\theta(t)) \cos{(\theta}(t))}_2^2 \text{ where } a(\theta), \theta'(\theta) \in V(\theta).
  \label{eq:minimization}
\end{equation}
Solution to this problem leads to the corresponding IMF, $s_1(t) = a(\theta(t)) \cos(\theta(t))$. 
To find the second and subsequent IMF's, we replace the signal $x(t)$ with the residual $r_k(t) = x(t) - \sum_{i = 1}^{k - 1} s_k(t)$.

To solve the NMP minimization problem an alternating scheme is used by fixing $\theta$ and minimizing over $a$, then fixing $a$ and updating $\theta$. It is also important to note that we first minimize 
\begin{equation}
  \argmin_{a, \theta} \abs{x(t) - a(\theta(t)) \cos{(\theta}(t))}_2^2 \text{ where } a(\theta), \theta'(\theta) \in V(\theta, \lambda),
  \label{eq:minlambda}
\end{equation}
with $\lambda = 0$ and slowly increase the value of $\lambda$ up to $1/2$. The corresponding pseudocode is in Algorithm \ref{alg:NMP}.

\begin{algorithm}[t]
\caption{Nonlinear Matching Pursuit (NMP) Method}
\begin{algorithmic}[1]
\STATE Input: measured signal $x(t)$ and initial guess for phase of IMF $\theta_0(t)$. For the 2nd and subsequent IMF's $x(t)$ is replaced with the remainder $r_k(t)$. 
\STATE Output: IMF $s(t)$

\STATE $\theta^0 := \theta_0(t)$, $\eta := 0$
\WHILE{$\eta < \lambda$}
\STATE $n = 0$
\WHILE{$n = 0$ or $\norm{\theta^{n+1} - \theta^n}_2 > \epsilon_0$}
\STATE{$a^{n+1},b^{n+1} := \arg\min_{a,b} \norm{x - a \cos(\theta^n) - b \sin(\theta^n)} _2^2$ s.t. $a(\theta),b(\theta) \in V(\theta,\eta)$} \COMMENT{Update $a$ and $b$}
\STATE{$\Delta \theta' := P_V(\theta;\eta) \left( \frac{d}{dt} \arctan \left( \frac{b(t)}{a(t)} \right) \right)$}
\STATE{$\Delta \theta(t) := \int_0^t \Delta \theta'(s) ds$}
\STATE{$\theta^{n+1} := \theta^n - \beta \Delta \theta$ where $\beta := \max \set{ \alpha \in [0,1]: \frac{d}{dt} \left( \theta^n - \alpha \Delta \theta \right) \geq 0 }$} \COMMENT{$\theta$ must be monotonic}
\STATE{$n := n + 1$}
\ENDWHILE
\STATE{$\eta = \eta + \Delta \eta$}
\ENDWHILE
\STATE{$a := \sqrt{\left(a^{n+1}\right)^2 + \left(b^{n+1}\right)^2}$}
\STATE{$\theta := \theta^{n+1}$}
\RETURN{$a \cos(\theta)$}
\end{algorithmic}{}
\label{alg:NMP}
\end{algorithm}


The minimization in \eqref{eq:minlambda} is nontrivial to compute. 
To solve it, Hou and Shi use the fact that projecting $a(\theta)$ into $V(\theta)$ is equivalent to applying a low pass filter in the $\theta$-coordinate \cite{hou2013data}. 
Pseudocode for this algorithm is shown in Algorithm~\ref{alg:NMP2}.
It is important to note that, in addition to taking in the measured signal $x(t)$, the minimization requires an initial guess for $\theta(t)$. 
Thus in the following sections, we will denote the NMP method as $\text{NMP}(x,\theta)$.

\subsection{The Blind Source Separation Problem} \label{subsec:factoring}
For systems with spatiotemporal dynamics, mixed time-series data from multiple sensors may be factored into a sum of spatiotemporal modes.
This problem is commonly known as blind source separation or factor analysis.

Generally, suppose we have spatiotemporal data $\bm{X} \in \R^{m \times n}$, which contains $n$ snapshots and $m$ measurement features at each snapshot. 
The goal is to decompose $\bm{X}$ into the product of two matrices $\bm{B} \in \R^{m \times r}$ and $\bm{S} \in \R^{r \times n}$ such that
\begin{equation}
\bm{X} = \bm{B} \bm{S}.
\label{eq:decompmat}
\end{equation} 
Equivalently, $\bm{X}$ may be expressed as 
\begin{equation}
  \bm{X} = \sum_{j = 1}^r b_j s_j,
\label{eq:decompsum}
\end{equation}
where $b_j \in \R^m$ is the $j^\text{th}$ column of  $\bm{B}$ and $s_j \in \R^n$ is the $j$th row of $S$. 
The column vectors $b_j$'s contain the spatial structures of the data, while the row vectors $s_j$'s contain the temporal structure. 
In other words, $b_j$'s are the spatial modes and $s_j$'s are the temporal modes of the data.
The dimension $r$ is typically chosen to optimize some objective; if $r < m$, then the decomposition can be used to reduce the dimensionality of $\bm{X}$ by representing the data in the basis of $b_j$'s.

This decomposition is highly under-determined; given different assumptions and constraints, a large variety of  different results for $\bm{B}$ and $\bm{S}$ may be obtained. 
For example, methods such as SVD and ICA make assumptions about the orthogonality and the statistical independence of the data, but they do not assume explicit temporal dynamics. 
Other methods such as DMD enforce a strict temporal structure on data with a linear dynamic model.  
Stronger assumptions restrict the types of data that can be models and reconstructed accurately; however, these assumptions, if appropriate for the system, may denoise the data and improve the interpretability of the results.
Sections \ref{subsubsec:SVD}, \ref{subsubsec:ICA}, and \ref{subsubsec:DMD} give an overview of these methods, and their properties are summarized in Table \ref{tab:algcomparison}.

\subsection{Singular Value Decomposition (SVD)} \label{subsec:svd}
\label{subsubsec:SVD}
One of the most widely used methods in matrix factorization is the SVD\footnote{Depending on the domain, this method (with small variations) is also known as Principal Component Analysis (PCA), Proper Mode Decomposition (POD), and the Karhunen-Lo$\tilde{\text{e}}$ve Decomposition, among others.}. 
Given a matrix $\bm{X} \in \R^{m \times n}$, SVD decomposes $\bm{X}$ into a product of three matrices
\begin{equation}
   \bm{X} = \bm{U} \bm{\Sigma} \bm{V}^T,
   \label{eq:svd}
 \end{equation} 
 where the left singular vectors $\bm{U} \in \R^{m \times m}$ and the right singular vectors $\bm{V} \in \R^{n \times n}$ are unitary matrices, and $\bm{\Sigma} \in \R^{m \times n}$ is diagonal \cite{golub1970singular,jolliffe1986principal}. 
 It is customary that the diagonal elements of $\bm{\Sigma}$ be expressed as
\[ \bm{\Sigma} = 
\begin{bmatrix}
    \sigma_1 & 0 & 0 & \dots  & 0 \\
    0 & \sigma_2 & 0 & \dots  & 0 \\
    \vdots & \vdots & \vdots & \ddots & \vdots \\
    0 & 0 & 0 & \dots  & \sigma_n \\
    0 & 0 & 0 & \dots & 0 \\
    \vdots & \vdots & \vdots & \vdots & \vdots \\
    0 & 0 & 0 & \dots & 0 \\
\end{bmatrix},
\]
where the singular values are in decreasing order, $\sigma_1 \geq \sigma_2 \geq \dots \sigma_n$. 
The rank of $\bm{X}$ is $R$, which corresponds to the number of nonzero $\sigma_i$'s \cite{kutz2013data}. Equivalently, we may incorporate the weightings $\bm{\Sigma}$ into $\bm{U}$. Letting $\bm{B} = \bm{U \Sigma}$ and $\bm{S} = \bm{V}^T$ we recover ~\eqref{eq:decompmat}. 


When $\bm{X}$ is spatiotemporal data, we may interpret $\bm{B}$ as the spatial modes of $\bm{X}$ and $\bm{S}$ as the temporal modes of $\bm{X}$.
Consider the matrix $\bm{X}_r$ defined as
\begin{equation}
  \bm{X}_r = \sum_{j = 1}^r b_j s_j
\end{equation}
where $0 \leq r \leq R$. 
This matrix has rank $r$;
importantly, this  $\bm{X}_r$ is the best rank $r$ approximation to $\bm{X}$. 
More precisely, if $\bm{Y}$ is a rank $r \leq R$ matrix, then $\| \bm{X} -  \bm{Y} \|$ is minimized for $\bm{Y} = \bm{X}_r$ with respect to both the $\ell_2$ and Frobenius norms. 
The relative error in the rank r approximation with respect to the $\ell_2$ norm is 
\begin{equation}
\frac{\| \bm{X} -  \bm{X}_r \|_2}{\norm{\bm{X}}_2} = \frac{\sigma_{r + 1}}{\sigma_1}.
\end{equation}
From this, we see that if the singular values $\sigma_i$ decay sufficiently rapidly such that $\sigma_{r+1} \ll \sigma_1$, then $\bm{X}_r$ will be a very good approximation to $\bm{X}$. 
This property makes SVD a popular tool for performing dimensionality reduction and mode extraction.

\subsection{Independent Component Analysis (ICA)} \label{subsec:ica}
\label{subsubsec:ICA}
Independent component analysis is another commonly used method for performing blind source separation. Common applications of ICA include brain imaging, finance, and image feature extraction \cite{back1997first,lee1999independent,hyvarinen2004independent,hyvarinen2013independent}. 
Like other blind signal separation problems, we assume that there is a set of $r$ signals $s_1, \ldots, s_r \in \R^{n}$ and we measure linear combinations $b_{i,j} \in \R$ to form the signals $x_i(t) \in \R^m$, as in~\eqref{eq:decompsum}
\begin{equation}
  x_i = \sum_{j = 1}^r b_{i,j} s_j.
\end{equation}

ICA makes the following assumptions:

\begin{enumerate}
  \item The source signals $s_j$ are mutually statistically independent. 
  \item The $s_j$'s follow non-gaussian distributions.
  \item The mixing matrix $\bm{B}$ is orthogonal.\\
\end{enumerate}

By the central limit theorem, a linear combination of signals tends to be more Gaussian than the distribution of a single signal. 
Consequently, many ICA algorithms, compute the $\bm{B}$ and $\bm{S}$ matrices by maximizing the nongaussianity of the source signals, which can be computed using kurtosis or negentropy. Alternatively, some methods optimize the statistical independence between signals by minimizing their mutual information or by employing joint diagonalization~\cite{cardoso1993blind,bell1995information,belouchrani1997blind}. 

One of the most popular algorithms for performing ICA is FastICA~\cite{hyvarinen2000independent}. 
To solve this problem, FastICA uses the method of projection pursuit: find the direction $w_j$ such that the projection $w_j^T \bm{X}$ maximizes the measure of nongaussianity. Constraining the $w_j$'s to be orthogonal yields a solution to $\bm{W} \bm{X} = \bm{S}$ (where the $j$th row of $\bm{W}$ is $w_j$). 
By construction, $\bm{W}$ is an orthogonal matrix, letting $\bm{B}$ be the pseudoinverse $\bm{W}^{\dagger} = \bm{W}^T$ to yield $\bm{X} = \bm{B} \bm{S}$ as desired.

\subsection{Dynamic Mode Decomposition (DMD)} \label{subsec:dmd}
\label{subsubsec:DMD}
As with SVD and ICA, we have some data $\bm{X} \in \R^{m \times n}$, but here we define $x(t_k) \in \R^m$ to be the state of the system at time $t_k$. 
We will further assume that the state has been sampled evenly in time at some spacing $\Delta t$ at a total of $n$ snapshots.
DMD has become a popular tool to model dynamical systems in the fields of fluid mechanics, neuroscience, and image analysis \cite{rowley2009spectral, schmid2010dynamic, kutz2016background, brunton2016extracting}.

The goal of DMD is to determine the best linear operator to $\bm{B} : \R^m \to \R^m$ such that
\begin{equation}
  x(t_{k + 1}) \approx  \bm{B} x(t_k).
\end{equation}
We let 
\[
 \bm{X}_1^{n - 1} = 
  \begin{bmatrix}
  | & | & \cdots & | \\
  x(t_1) & x_2(t_2) & \cdots & x(t_{n-1}) \\
  | & | & \cdots & | \\
\end{bmatrix} \text{ , }
\bm{X}_2^n = 
  \begin{bmatrix}
  | & | & \cdots & | \\
  x(t_2) & x(t_3) & \cdots & x(t_n) \\
  | & | & \cdots & | \\
\end{bmatrix}.
\]
Then, we can equivalently define $\bm{B} \in \R^{m \times m}$ to be the operator such that
\begin{equation}
	\bm{X}_2^n \approx \bm{B} \bm{X}_1^{n - 1}.
\end{equation} 
To find $\bm{B}$, we must then solve the minimization problem
\begin{equation}
  \min_{\bm{B}} \norm{\bm{X}_2^n - \bm{B} \bm{X}_1^{n-1}}_F,
\end{equation}
where $\norm{\cdot}_F$ denotes the Frobenius norm. 
A unique solution to this problem can be obtained using the exact DMD method~\cite{Tu2014jcd}. 
For noisy data, we can use more robust methods such as optimized DMD~\cite{askham2018variable}.

One key benefit of DMD is that it builds an explicit temporal model, which allows future state prediction. Specifically, if we let $\left\{ \lambda_j \right\}$ and $\left\{ b_j \right\}$ be the eigenvalues and eigenvectors of $\bm{B}$, respectively, then 
\begin{equation}
	x(t) = \sum_{j = 1}^r b_j e^{\omega_j t} c_j
	\label{eq:DMD1}
\end{equation}
where $\omega_j = \ln(\lambda_j) / \Delta t$ and $c \in \R^m$ corresponds to the initial conditions of the state. Thus, to compute the state at an arbitrary time $t$, simply evaluate \eqref{eq:DMD1} at that time.

Defining $s_j(t) = \exp(\omega_j t) c_j$, then we have
\begin{equation}
  x(t_k) = \sum_{j = 1}^r b_j s_j (t_k).
\end{equation}
Thus, we can think of DMD as a matrix factorization as in~\eqref{eq:decompsum}, where the spatial modes are the $b_j$'s, and the temporal modes $s_j$ all have complex exponential temporal dependence. 

\section{Spatiotemporal Intrinsic Mode Decomposition (STIMD)}
\label{sec:method}

Our STIMD algorithm leverages the mathematical and algorithmic structures of EMD and NMP for improved spatiotemporal decompositions.  
The mathematical framework and algorithmic implementation are given in the following subsections.
All the code developed to implement STIMD, along with scripts to reproduce results in the figures, are openly available at \url{https://github.com/BruntonUWBio/STIMD}.

\subsection{Method Description}
As introduced in Section~\ref{sec:relatedwork}, we assume we have a set of source signals $s_i$ that are linearly mixed to form the observed signals $x_i$,
\begin{equation}
  x_i(t) = \sum_{j = 1}^r b_{i,j} s_j(t),
\end{equation}
or equivalently in matrix form 
\begin{equation}
  \bm{X} = \bm{B} \bm{S}.
\end{equation}

Here we assume that all modes $s_j$ are IMF's with interwave frequency modulation, as defined in Section~\ref{subsec:nmp}. 
Thus, each $s_j$ takes the form
\begin{equation}
  s_j(t) = a_j(\theta_j(t)) \cos(\theta_j(t)) = a_j(\theta_j) \cos(\theta_j),
\end{equation}
where $a_j(\theta_j), \theta_j' \in V(\theta_j)$ are as defined in (\ref{eq:Vtheta}).

As summarized in Table~\ref{tab:algcomparison}, this mathematical architecture provides a compromise between ICA, SVD and DMD. 
We obtain a model of the temporal dynamics that are constrained to obey a set of dynamics commonly found in physical systems, but they are not restricted to be stationary signals as in DMD.
In addition, we have the ability to compute a Hilbert spectrum and perform future state prediction of non-stationary signals.

\begin{table}[]
\centering
\caption{Comparison of assumptions, features, and limitations of the spatiotemporal decomposition algorithms. }
\label{tab:algcomparison}
\begin{tabular}{@{}lcccc@{}} 
\toprule
                         & SVD        & ICA        & DMD        &  STIMD  \\
\midrule                         
Derives orthogonal spatial modes & \checkmark & \checkmark & \xmark     & \xmark  \\ 
Models temporal dynamics  & \xmark     & \xmark     & \checkmark & \checkmark   \\ 
Adapts to non-stationary signals   & \checkmark & \checkmark & \xmark     & \checkmark \\ 
Predicts future states  & \xmark     & \xmark     & \checkmark & \checkmark   \\ 
Suitable for Hilbert spectrum computation & \xmark     & \xmark     & \checkmark & \checkmark \\ 
\bottomrule
\end{tabular}
\end{table}

\subsection{Method Implementation}
Inspired by the FastICA algorithm~\cite{hyvarinen2000independent}, we propose a method based on projection pursuit. Our goal is to find the direction $b_1$ such that $b_1^T X$ is an IMF $s_1$. If we apply the NMP algorithm (see Secion~\ref{app:nmp}) to an IMF, we expect NMP to return that same IMF. 
Thus, we seek the fixed point
\begin{equation}
  s_1^T = \text{NMP}(b_1^T \bm{X}, \theta_1) \approx b_1^T \bm{X}.
\end{equation}
Without loss of generality, we constrain $b_1$ to have unit norm. 
Applying the pseudoinverse of $b_1$ to both sides gives
\begin{equation}
 b_1 \text{NMP}(b_1^T \bm{X}, \theta_1 ) \approx \bm{X}.
\end{equation}
To find $b_1$, we minimize the difference between the left and right hand sides
\begin{equation}
  b_1 = \argmin_{w} \norm{w \text{NMP}(w^T \bm{X},\theta_1) - \bm{X}}_2^2 \text{~~~subject to $\norm{w}_2 = 1$}.
\end{equation}
To find the next IMF, we apply the same method to the remainder $\bm{R}_1$
\begin{equation}
  \bm{R}_1 = \bm{X} - b_1 s_1^T.
\end{equation}
Thus, to find the $i$th IMF $s_i$ we solve the minimization problem 
\begin{equation}
  b_i = \argmin_{w} \norm{w \text{NMP}(w^T \bm{R}_{i-1},\theta_1) - \bm{R}_{i-1}}_2^2 \text{~~~subject to $\norm{w}_2 = 1$},
\end{equation}
where $\bm{R}_i$ is given by
\begin{equation}
  \bm{R}_i = \bm{X} - \sum_{j = 1}^i b_j s_j^T .
\end{equation}
The $i$th IMF $s_i$ is then
\begin{equation}
  s_i = \text{NMP}(b_i^T \bm{R}_{i-1}).
\end{equation}
This decomposition in matrix form produces
\begin{equation}
  \bunderbrace{\begin{bmatrix}
  - & x_1 & -  \\
  - & x_2 & - \\
  \vdots & \vdots & \vdots \\
  - & x_m & -
  \end{bmatrix}}{\bm{X}} \approx
  \bunderbrace{\begin{bmatrix}
  | & | & \cdots & | \\
  b_1 & b_2 & \cdots & b_r \\
  | & | & \cdots & | 
  \end{bmatrix}}{\bm{B}}
  \bunderbrace{\begin{bmatrix}  - & s_1 & -  \\
  - & s_2 & - \\
  \vdots & \vdots & \vdots \\
  - & s_r & -
  \end{bmatrix}}{\bm{S}}
  \label{eq:myalg}.
\end{equation}
The corresponding pseudocode for STIMD is shown in Algorithm~\ref{alg:stimd}. 
An implementation in Python can be found at \url{https://github.com/BruntonUWBio/STIMD}.

\begin{algorithm}
\caption{STIMD}
\begin{algorithmic}
\STATE Input: measured signals $\bm{X} = [x_1, x_2, \ldots, x_n]^T$ and initial guess for phases of IMF's $\bm{\theta} = [\theta_1, \theta_2, \ldots, \theta_n]^T$
\STATE Output: Matrix of IMF's $\bm{S} = [s_1, s_2, \ldots, s_n]^T$ and mixing matrix  $\bm{B}$.
\STATE{$\bm{R}_0 := \bm{X}$}
\FOR{$i \in \left\{1,\ldots,n \right\}$}
\STATE{$b_i := \arg \min_w \norm{w \text{NMP}(w^T \bm{R}_{i-1},\theta_i) - \bm{R}_{i-1}}_2^2$ s.t. $\norm{w}_2 = 1$ }
\STATE{$s_i := \text{NMP}(b_i^T \bm{R}_{i-1}) $}
\STATE $\bm{R}_i := \bm{R}_{i - 1} -  b_i s_i^T$
\ENDFOR
\STATE $\bm{B} := [b_1, b_2, \ldots, b_n]$ 
\STATE $\bm{S} := [s_1^T, s_2^T, \ldots, s_n^T]^T$
\RETURN $\bm{B},\bm{S}$
\end{algorithmic}
\label{alg:stimd}
\end{algorithm}




\subsubsection{Future State Prediction}
\label{subsubsec:futurestate}
As a consequence of the connection to NMP, it is easy to extract $\theta_i$ evaluated in time from STIMD. 
From (\ref{eq:Vtheta}), we know that for each IMF $s_i(\theta) = a(\theta) \cos(\theta)$ we have assumed that $\theta'(\theta) \in V(\theta)$. 
In particular,
\begin{equation}
  \theta'(\theta) = \alpha_0 + \sum_k \beta_k \cos\left(\frac{k \theta}{L_\theta} \right) + \gamma_k \sin\left(\frac{k \theta}{L_\theta} \right).
\end{equation}
The coefficients $\alpha_0, \beta_k$, and $\gamma_k$ can all be computed using the Fourier transform in the $\theta$ coordinate. 
With these coefficients, we now have an implicit first order differential equation that can be integrated to compute $\theta(t)$, forming a model of $s_i(\theta(t))$ at any future time.


\subsection{Experiments on Synthetic Examples}

We demonstrate STIMD on a number of synthetic examples.
In addition to showing its ability to accurately extract mixed non-stationary source signals, we characterize its sensitivity to noise and initial conditions.


\begin{figure}
\centering
\begin{overpic}[scale=0.35,unit=1bp]{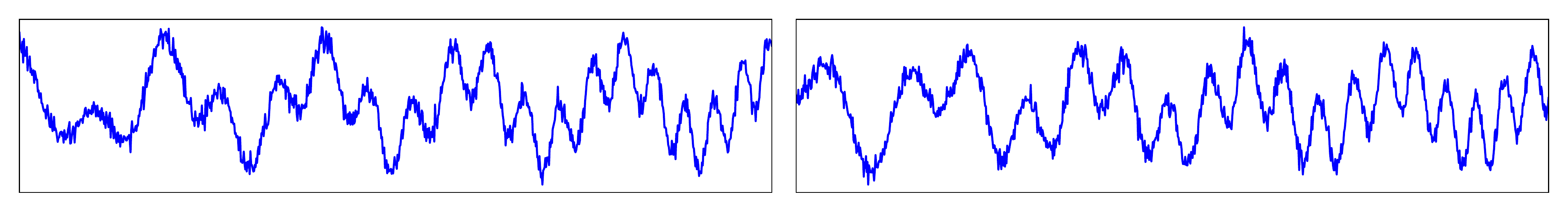}
\put(-6,1){\rotatebox{90}{\Centerstack{Observed \\ Signals}}}
\end{overpic}
\par\medskip 
\begin{overpic}[scale=0.35,unit=1bp]{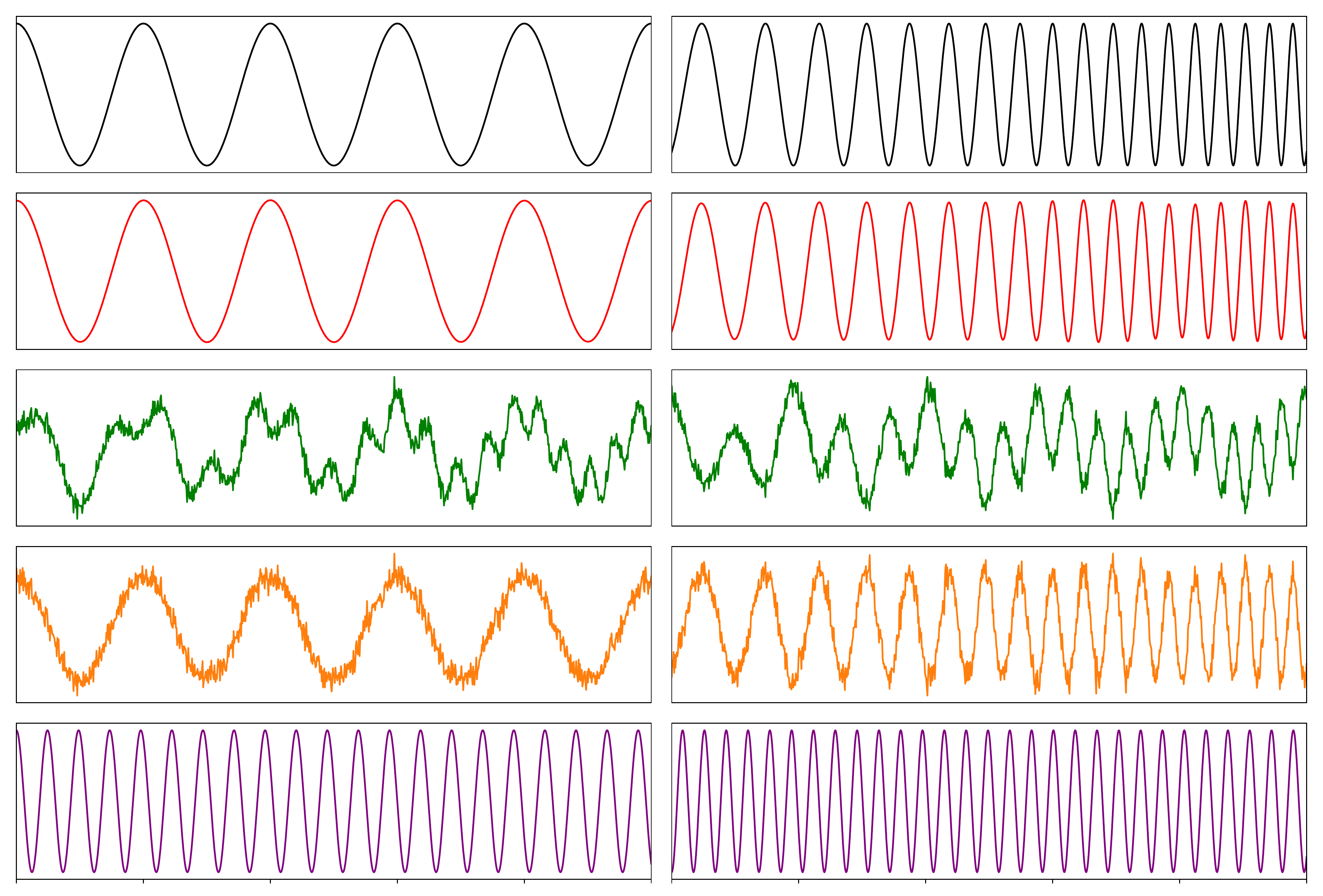}
\put(20,67){Mode 1}
\put(70,67){Mode 2}
\put(-6,56){\rotatebox{90}{\parbox{.5in}{\Centerstack{Source \\ Signals}}}}
\put(-2,43){\rotatebox{90}{STIMD}}
\put(-2,32){\rotatebox{90}{SVD}}
\put(-2,18){\rotatebox{90}{ICA}}
\put(-2,4){\rotatebox{90}{DMD}}
\put(20,-4){Time [s]}
\put(70,-4){Time [s]}

\put(0,-1){\small 0.0}
\put(9,-1){\small 0.1}
\put(19,-1){\small 0.2}
\put(28,-1){\small 0.3}
\put(38,-1){\small 0.4}
\put(46.5,-1){\small 0.5}

\put(50.5,-1){\small 0.0}
\put(59,-1){\small 0.1}
\put(68,-1){\small 0.2}
\put(78,-1){\small 0.3}
\put(87.5,-1){\small 0.4}
\put(96.5,-1){\small 0.5}
\end{overpic}
\vspace*{5mm}
\caption{Comparison of temporal modes extracted by several blind signal separation algorithms for a two-dimensional system. 
Top row: The observed mixed measurement signals in blue. These correspond to linear combinations of source signals plus a small amount of Gaussian-distributed noise.
Row 2: The two true source signals in black.   
Row 3: Modes extracted by STIMD in red. 
Row 4--6: Modes extracted by ICA (green), SVD (orange), and by optimized DMD (purple), respectively.}
  \label{fig:2dcomparison}
\end{figure}
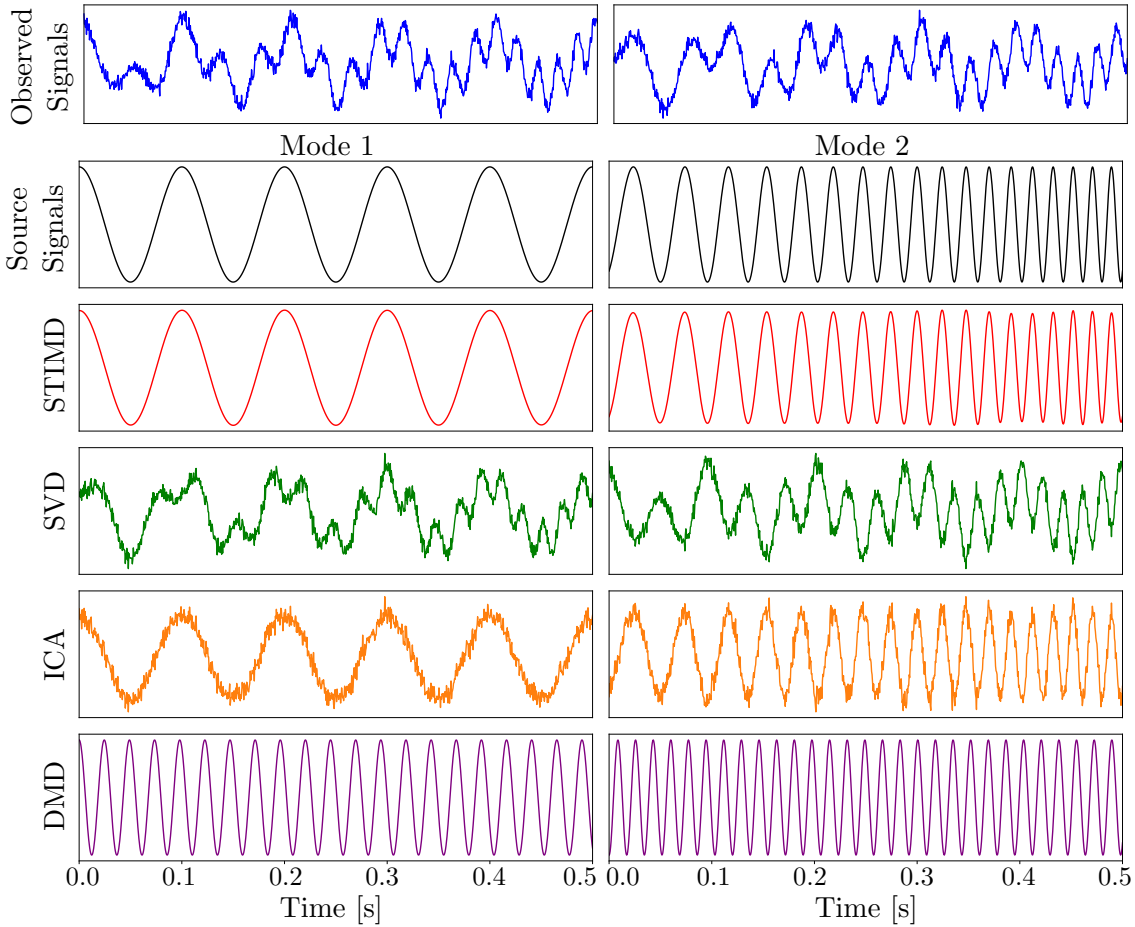

\subsubsection{2D Example}
\label{subsubsec:2d}
First, we consider a simple two-dimensional example. 
The true signals $s_1$ and $s_2$ are oscillatory signals with frequency modulation,
\begin{align*}
  s_1(t) =& \sin(10 \pi t) \\
  s_2(t) =& \sin(20 \pi (t+0.4)^2),
\end{align*}
for $t \in [0,1]$, and the mixing matrix $\bm{B}$ is given by 
\[
\bm{B} = \begin{bmatrix} 
            \cos(\phi_0) & -\sin(\phi_0) \\
            \sin(\phi_0) & \cos(\phi_0) \\
            \end{bmatrix},
\]
where $\phi_0 = 0.7$. 
By construction, $\bm{B}$ is orthogonal, since the matrix used is a standard rotation matrix. 
Adding a small amount of measurement noise to each of the sensors $N_1(t), N_2(t) \sim \mathcal{N}(0,0.1)$, the observed (mixed signals) $\bm{X} = \bm{B S} + \bm{N}$ are given by
\begin{align*}
x_1(t) =& \cos(\phi) s_1(t) -\sin(\phi) s_2(t) + N_1(t) \\
x_2(t) =& \sin(\phi) s_1(t) +\cos(\phi) s_2(t) + N_2(t).
\end{align*}


Figure~\ref{fig:2dcomparison} shows the source signals ($s_1, s_2$) (black), mixed noisy measured signals $x_1(t), x_2(t)$ (blue) and the signals extracted using STIMD (red). 
For comparison, the modes extracted from SVD (green), ICA (orange), and optimized DMD (purple) are also showh. 
From this, we see that the STIMD modes closely capture the true source signals. 
The SVD modes are clearly still a mixture of the measured signals. 
The ICA modes contain some amplitude modulation in time not seen in the true signals. 
More importantly, the modes extracted by SVD and ICA are not IMF's and consequently are not guaranteed to have meaningful Hilbert spectrums. 
Lastly, using the optimized DMD algorithm, neither non-stationary mode is extracted correctly.

\subsubsection{3D Example}

Next, we consider an example containing 3 modes. The source signals are
\begin{align*}
s_1(t) =& \cos(20 \pi t -5 \sin(\pi t)) \\
s_2(t) =& \cos(60 \pi t+2 \sin(4 \pi t)) \\
s_3(t) =& \cos(90 \pi t+3 \sin(8 \pi t)) , 
\end{align*} 
and the mixing matrix is
\[
  \bm{B} = \begin{bmatrix}
              \cos(\phi_1) \sin(\phi_2) & -\sin(\phi_1) & \cos(\phi_1) \cos(\phi_2) \\
              \sin(\phi_1) \sin(\phi_2) & \cos(\phi_1) &\sin(\phi_1) \cos(\phi_2)\\
              \cos(\phi_2) & 0 & -\sin(\phi_2)
      \end{bmatrix},
\]
with $\phi_1 = 0.6$ and $\phi_2 = 0.7$. 
As before, the observed signals $x_i$ are the linear combinations 
\[
\underbrace{\begin{bmatrix}
- & x_1(t) & - \\
- & x_2(t) & - \\
- & x_3(t) & - \\
\end{bmatrix}}_{\text{\normalsize $\bm{X}$}}
= 
\bm{B}
\underbrace{\begin{bmatrix}
- & s_1(t) & - \\
- & s_2(t) & - \\
- & s_3(t) & - \\
\end{bmatrix}}_{\text{\normalsize $\bm{S}$}},
\]
or writing it more succinctly, $\bm{X} = \bm{B} \bm{S}$.

\begin{figure}
\centering
\begin{overpic}[scale=0.3,unit=1bp]{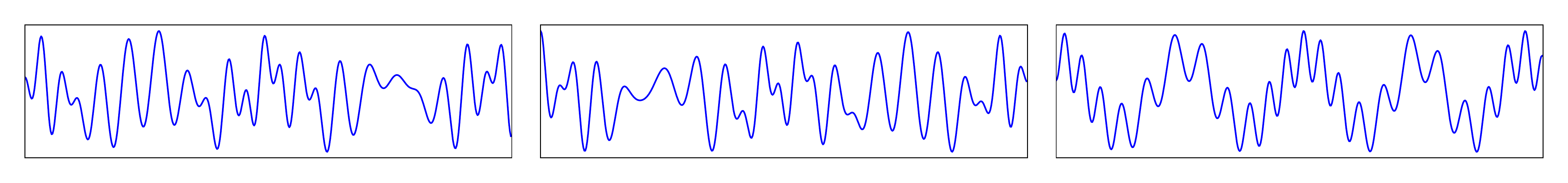}
\put(10,62){(a)}
\put(60,62){(b)}
\put(10,35){(c)}
\put(60,35){(d)}
\put(-6,1){\rotatebox{90}{\Centerstack{Observed \\ Signals}}}
\end{overpic}
\par\medskip 
\begin{overpic}[scale=0.3,unit=1bp]{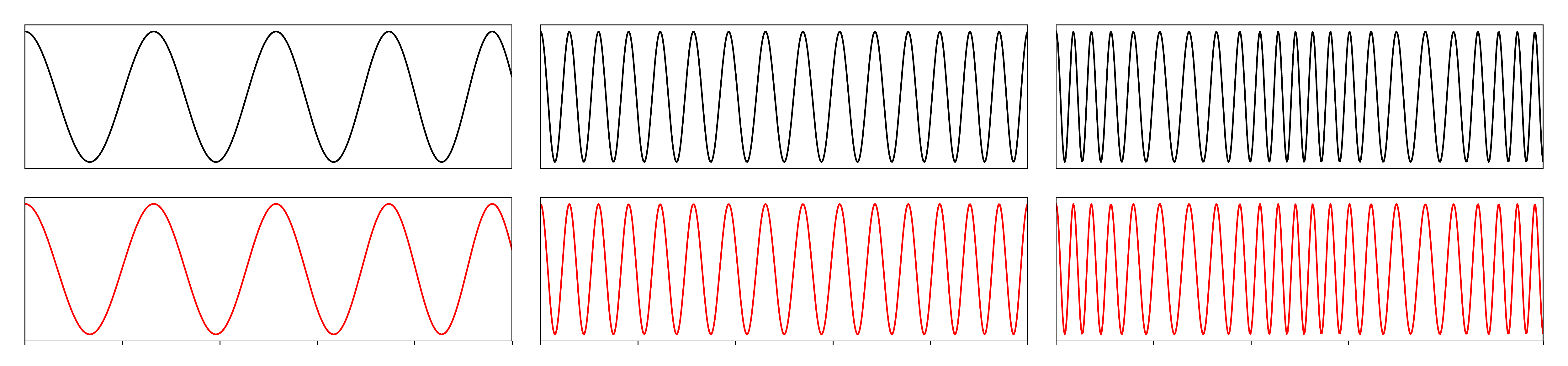}
\put(15,23){Mode 1}
\put(46,23){Mode 2}
\put(80,23){Mode 3}
\put(-6,14){\rotatebox{90}{\Centerstack{Source \\ Signals}}}
\put(-2,2){\rotatebox{90}{STIMD}}
\put(15,-4.){Time [s]}
\put(46,-4.){Time [s]}
\put(80,-4. ){Time [s]}
\put(1,-1){\small 0.0}
\put(7,-1){\small 0.1}
\put(12,-1){\small 0.2}
\put(18,-1){\small 0.3}
\put(24,-1){\small 0.4}
\put(30,-1){\small 0.5}

\put(34,-1){\small 0.0}
\put(40,-1){\small 0.1}
\put(46,-1){\small 0.2}
\put(51,-1){\small 0.3}
\put(57,-1){\small 0.4}
\put(63,-1){\small 0.5}

\put(67,-1.){\small 0.0}
\put(73,-1.){\small 0.1}
\put(79,-1.){\small 0.2}
\put(85,-1.){\small 0.3}
\put(90,-1.){\small 0.4}
\put(96,-1.){\small 0.5}
\end{overpic}
\vspace*{5mm}
\caption{Example of STIMD applied to spatiotemporal data with three modes. 
Top: Observed signals which are  linear combinations of source signals. 
Center: The source signals. 
Bottom: Signals reconstructed using STIMD.}
\label{fig:test1}
\end{figure}

The original modes ($s_1,s_2,s_3$) the observed signals ($x_1,x_2,x_3$) and the modes extracted by STIMD are shown in Figure \ref{fig:test1}. The modes extracted by STIMD are nearly identical to the original signals.

\begin{figure}[th]
\centering
\begin{overpic}[scale=0.3,unit=1bp]{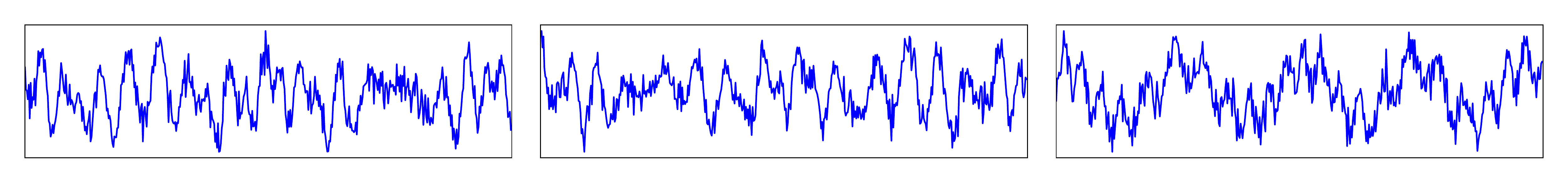}
\put(-6,1){\rotatebox{90}{\Centerstack{Observed \\ Signals}}}
\end{overpic}
\par\medskip 
\begin{overpic}[scale=0.3,unit=1bp]{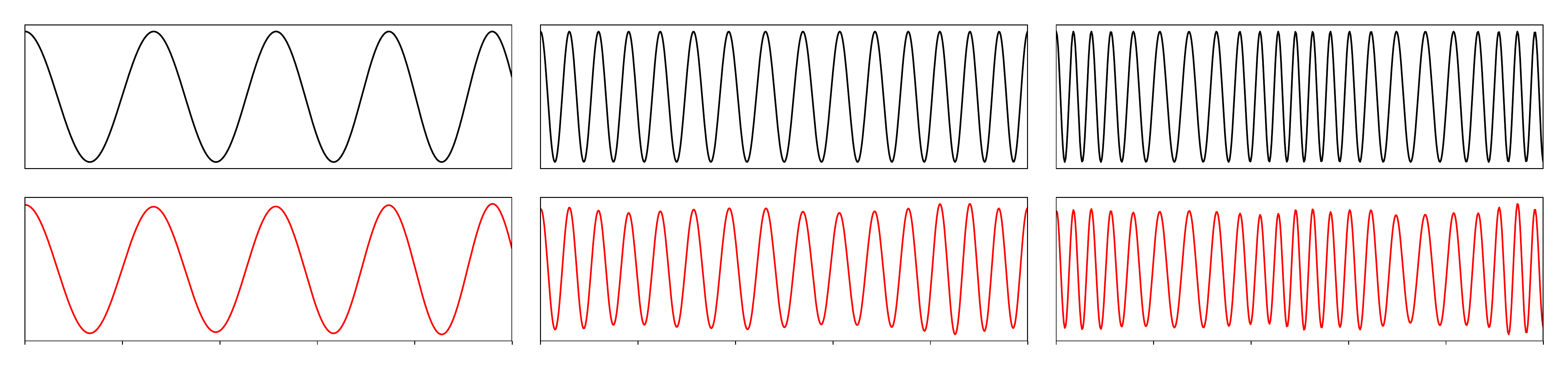}
\put(15,23){Mode 1}
\put(46,23){Mode 2}
\put(80,23){Mode 3}
\put(-6,14){\rotatebox{90}{\Centerstack{Source \\ Signals}}}
\put(-2,2){\rotatebox{90}{STIMD}}
\put(15,-4.){Time [s]}
\put(46,-4.){Time [s]}
\put(80,-4. ){Time [s]}
\put(1,-1){0.0}
\put(7,-1){\small 0.1}
\put(12,-1){\small 0.2}
\put(18,-1){\small 0.3}
\put(24,-1){\small 0.4}
\put(30,-1){\small 0.5}

\put(34,-1){\small 0.0}
\put(40,-1){\small 0.1}
\put(46,-1){\small 0.2}
\put(51,-1){\small 0.3}
\put(57,-1){\small 0.4}
\put(63,-1){\small 0.5}

\put(67,-1.){\small 0.0}
\put(73,-1.){\small 0.1}
\put(79,-1.){\small 0.2}
\put(85,-1.){\small 0.3}
\put(90,-1.){\small 0.4}
\put(96,-1.){\small 0.5}
\end{overpic}
\vspace{5mm}
\caption{Example of STIMD applied to spatiotemporal data with three modes in the case of measurement noise. Top: The source signals. Center: Observed signals which are  linear combinations of source signals plus Gaussian distributed noise. Bottom: Signals reconstructed using STIMD.}
\label{fig:test2}
\end{figure}

\begin{figure}[tb]
  \centering
  \includegraphics[scale=0.45]{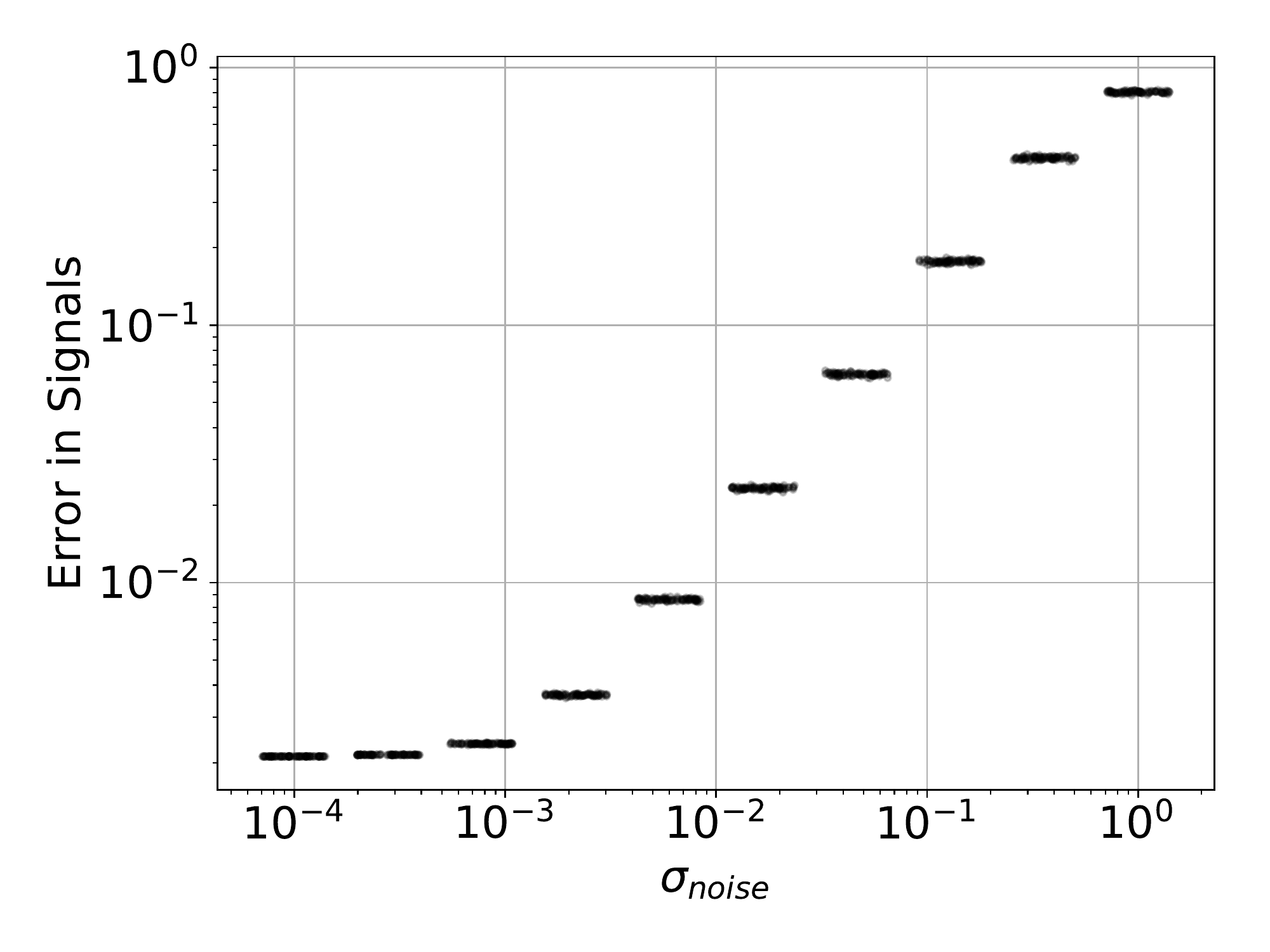}
  \caption{Characterization of accuracy of STIMD modes as a function of measurement noise. For each given value of noise 100 realizations of Gaussian distributed noise (with standard deviation $\sigma$) were added to the true signal. The relative error between the IMF's extracted by STIMD and the true IMF's are plotted. Random uniform jitter is added on a log-scale to each set of trials for visualization.}
  \label{fig:ErrorInSignals}
\end{figure}

We next consider how the results are affected by noise.
As an example, we add Gaussian noise with standard deviation $\sigma = 0.3$ to the measurements like for the 2D example. 
As in the noiseless case, we plot the original signals $s_i$ the measured signals $x_i$ and the STIMD modes (Figure~\ref{fig:test2}). 
There is a small amount of amplitude modulation not present in the original signals (which all have amplitude 1). 
Even so, the frequencies are nearly identical to the true signals. 

Figure~\ref{fig:ErrorInSignals} characterizes how the results are affected for various noise levels.  
Here we add Gaussian noise to the measured signals, sweeping standard deviations $\sigma$ from $10^{-4}$ to $1$. 
For each value of $\sigma$ we perform 100 trials and record the relative error of the extracted modes to the source signals with respect to the 2-norm. 
As expected, the error increases for increasing $\sigma$ and only becomes order $1$ (about the same size as the signals) when $\sigma$ is of the same order of magnitude.  

\begin{figure}[h]
  \centering
  \vspace{5mm}
  \begin{overpic}[scale=0.3,unit=1bp]{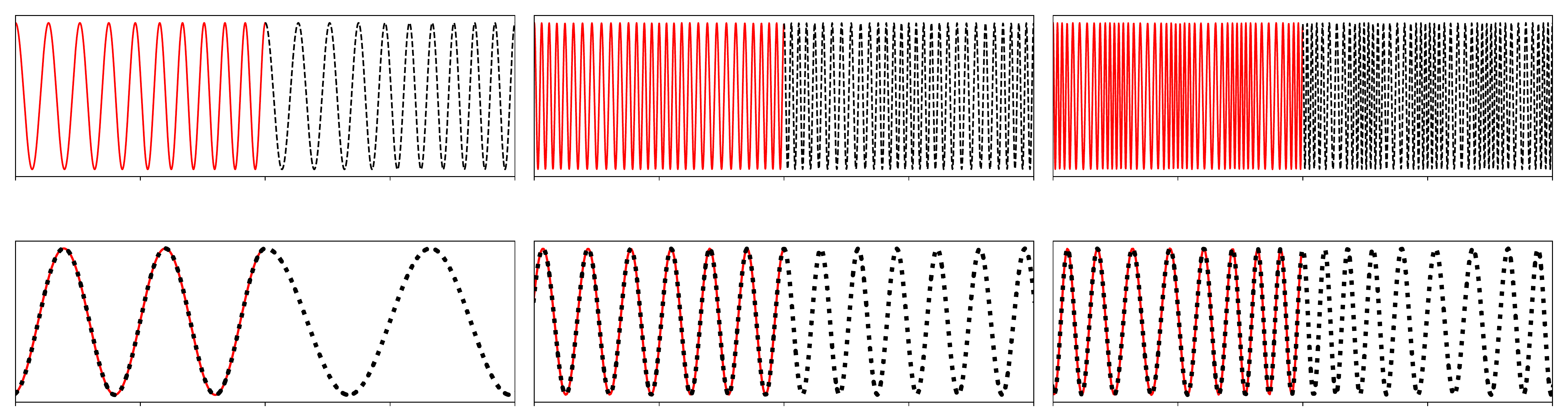}
\put(15,26.5){Mode 1}
\put(46,26.5){Mode 2}
\put(80,26.5){Mode 3}

\put(0.5,13.5){\small 0.0}
\put(7.5,13.5){\small 0.5}
\put(15,13.5){\small 1.0}
\put(23.5,13.5){\small 1.5}
\put(30,13.5){\small 2.0}

\put(34,13.5){\small 0.0}
\put(40.5,13.5){\small 0.5}
\put(48.5,13.5){\small 1.0}
\put(56.5,13.5){\small 1.5}
\put(63,13.5){\small 2.0}

\put(67,13.5){\small 0.0}
\put(73.5,13.5){\small 0.5}
\put(81.5,13.5){\small 1.0}
\put(89.5,13.5){\small 1.5}
\put(96.5,13.5){\small 2.0}

\put(15,-4.){Time [s]}
\put(46,-4.){Time [s]}
\put(80,-4. ){Time [s]}

\put(0.5,-1){\small 0.8}
\put(7.5,-1){\small 0.9}
\put(15,-1){\small 1.0}
\put(23.5,-1){\small 1.1}
\put(30,-1){\small 1.2}

\put(34,-1){\small 0.8}
\put(40.5,-1){\small 0.9}
\put(48.5,-1){\small 1.0}
\put(56.5,-1){\small 1.1}
\put(63,-1){\small 1.2}

\put(67,-1.){\small 0.8}
\put(73.5,-1.){\small 0.9}
\put(81.5,-1.){\small 1.0}
\put(89.5,-1.){\small 1.1}
\put(96.5,-1.){\small 1.2}

\put(15,-4.){Time [s]}
\put(46,-4.){Time [s]}
\put(80,-4. ){Time [s]}
  \end{overpic}
  \vspace*{5mm}
  \caption{Example of future state prediction for spatiotemporal data with three modes. Observed signals, and corresponding STIMD modes are shown in Figure \ref{fig:test1}. 
  Top row: The STIMD modes (red) were computed over the interval $t \in [0,1]$, while the state was predicted over the window $t \in [0,2]$ (black dotted). Note that the prediction ranges over the full interval. However, for illustration only the future state is shown.
  Bottow row: Same system as in the top row, zoomed in over the interval $t \in [0.8, 1.2]$.}
  \label{fig:futurestate}
\end{figure}

Figure~\ref{fig:futurestate} illustrates how the STIMD results from this three-dimensional system can be used for future state prediction, as described in Section~\ref{subsubsec:futurestate}.
As an example, consider the STIMD modes (red) over the range $t \in [0,1]$ for the noiseless system(see Figure \ref{fig:test1}).
We compute and predict the modes (black dotted lines) over the greater range of $t \in [0, 2]$ (Figure \ref{fig:futurestate}). 
The predicted modes follow the STIMD modes accurately over the interval $[0,1]$ and accounts accurately for the frequency modulation and amplitude modulation over $[1,2]$. 
Note that the state prediction will be only as good as the original reconstruction. 
In the case of noise the error in the prediction will go as the error in the STIMD modes.

\subsection{4D Example}

As a final synthetic example, we show STIMD applied to the case of eight (8) observed signals and four (4)source signals
\begin{align*}
  s_1(t) =& \cos(20 \pi t - 5 \sin(\pi t)) \\
  s_2(t) =& \cos(30 \pi t + \sin(4 \pi t)) \\
  s_3(t) =& \cos(60 \pi t + 3 \sin(5 \pi t)) \\
  s_4(t) =& \cos(80 \pi t + 4 \sin(5 \pi t)).
\end{align*}
Thus, the mixing matrix $\bm{B} \in \R^{8 \times 4}$ is a non-square matrix. 
Noise with standard deviation $\sigma = 0.1$ is added. 

\begin{figure}[htbp]
\centering
\begin{overpic}[scale=0.3,unit=1bp]{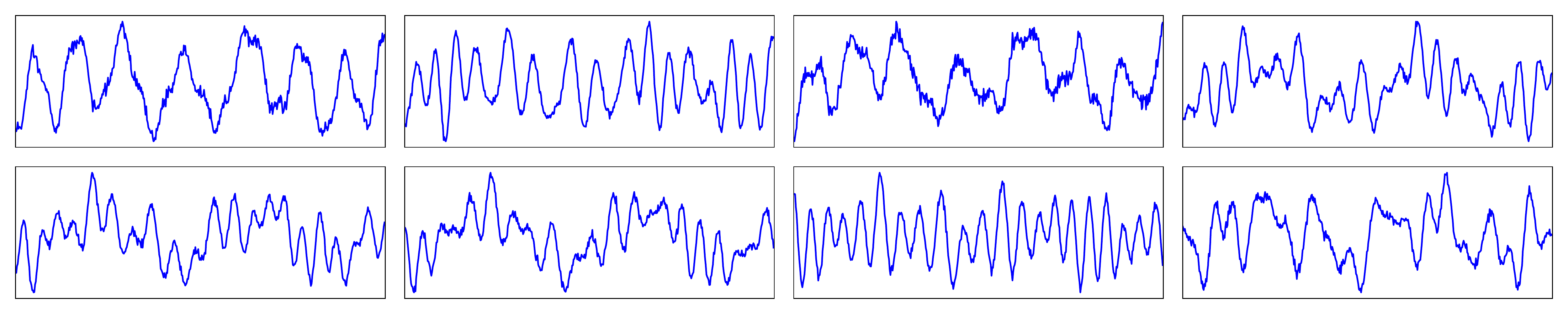}
\put(-6,5){\rotatebox{90}{\Centerstack{Observed \\ Signals}}}
\end{overpic}
\par\medskip 
\begin{overpic}[scale=0.3,unit=1bp]{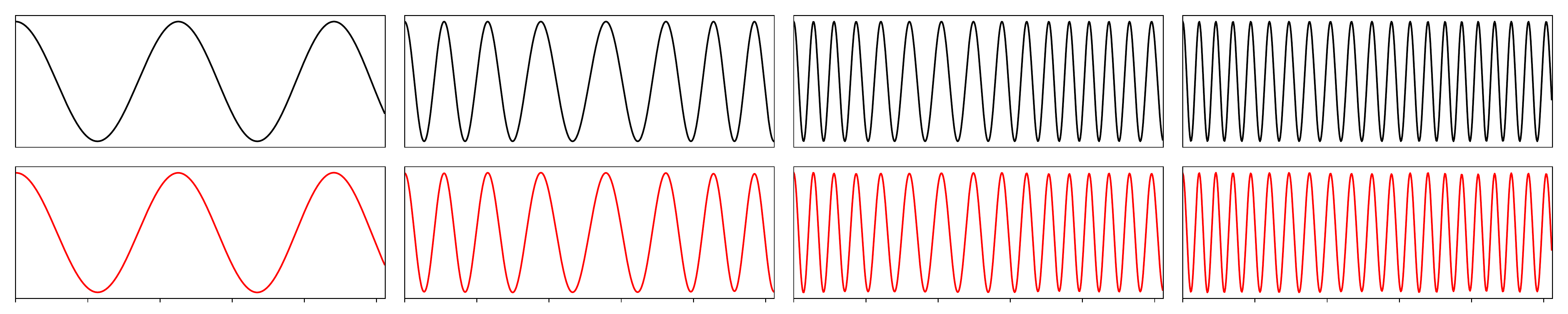}
\put(8,20){Mode 1}
\put(34,20){Mode 2}
\put(58,20){Mode 3}
\put(82,20){Mode 4}
\put(-6,12){\rotatebox{90}{\Centerstack{Source \\ Signals}}}
\put(-2,1){\rotatebox{90}{STIMD}}

\put(0,-1){\footnotesize 0.0}
\put(4,-1){\footnotesize 0.1}
\put(9,-1){\footnotesize 0.2}
\put(13.5,-1){\footnotesize 0.3}
\put(18,-1){\footnotesize 0.4}
\put(22,-1){\footnotesize 0.5}

\put(25,-1){\footnotesize 0.0}
\put(29,-1){\footnotesize 0.1}
\put(33.5,-1){\footnotesize 0.2}
\put(38,-1){\footnotesize 0.3}
\put(43,-1){\footnotesize 0.4}
\put(47,-1){\footnotesize 0.5}

\put(50,-1){\footnotesize 0.0}
\put(54,-1){\footnotesize 0.1}
\put(58.5,-1){\footnotesize 0.2}
\put(63,-1){\footnotesize 0.3}
\put(68,-1){\footnotesize 0.4}
\put(72,-1){\footnotesize 0.5}

\put(75,-1){\footnotesize 0.0}
\put(79,-1){\footnotesize 0.1}
\put(83.5,-1){\footnotesize 0.2}
\put(88,-1){\footnotesize 0.3}
\put(93,-1){\footnotesize 0.4}
\put(97,-1){\footnotesize 0.5}

\put(8,-4.){Time [s]}
\put(34,-4.){Time [s]}
\put(58,-4. ){Time [s]}
\put(82,-4. ){Time [s]}

\end{overpic}
\vspace*{7mm}
  \caption{Example of STIMD applied to spatiotemporal data with four modes in the case of measurement noise, and eight measured signals. Top row: The source signals. 2nd row and 3rd row: Observed signals which are  linear combinations of source signals plus Gaussian distributed noise. Bottom row: Signals reconstructed using algorithm.}
  \label{fig:Test4D}
\end{figure}

The source signals $\bm{S}$, the observed signals $\bm{X}$, and the STIMD modes are shown in Figure~ \ref{fig:Test4D}. 
Besides a small amount of amplitude modulation, once again the STIMD modes are nearly identical to the original signals.

\subsection{Initial Conditions}
It is important to emphasize that the NMP algorithm, and consequently the STIMD algorithm, takes as input initial guesses for the phases $\theta_i(t)$ of the IMF's.  
In many cases, only coarse guesses are needed. 
For example, for the 3D example, guesses corresponding to the central frequencies, $\theta_1(t) = 20 \pi t$ and $\theta_2(t) = 60 \pi t$ and $\theta_3(t) = 90 \pi t$ are needed. 
For one dimension, Shi recommends taking the Fourier transform and picking the peaks in the spectrum as initial guesses~\cite{hou2013data}. 
For the STIMD algorithm, we recommend using the peaks in the Fourier spectra of the temporal modes computed using FasTICA. 

\begin{figure}[t]
  \centering
  \includegraphics[scale=0.40]{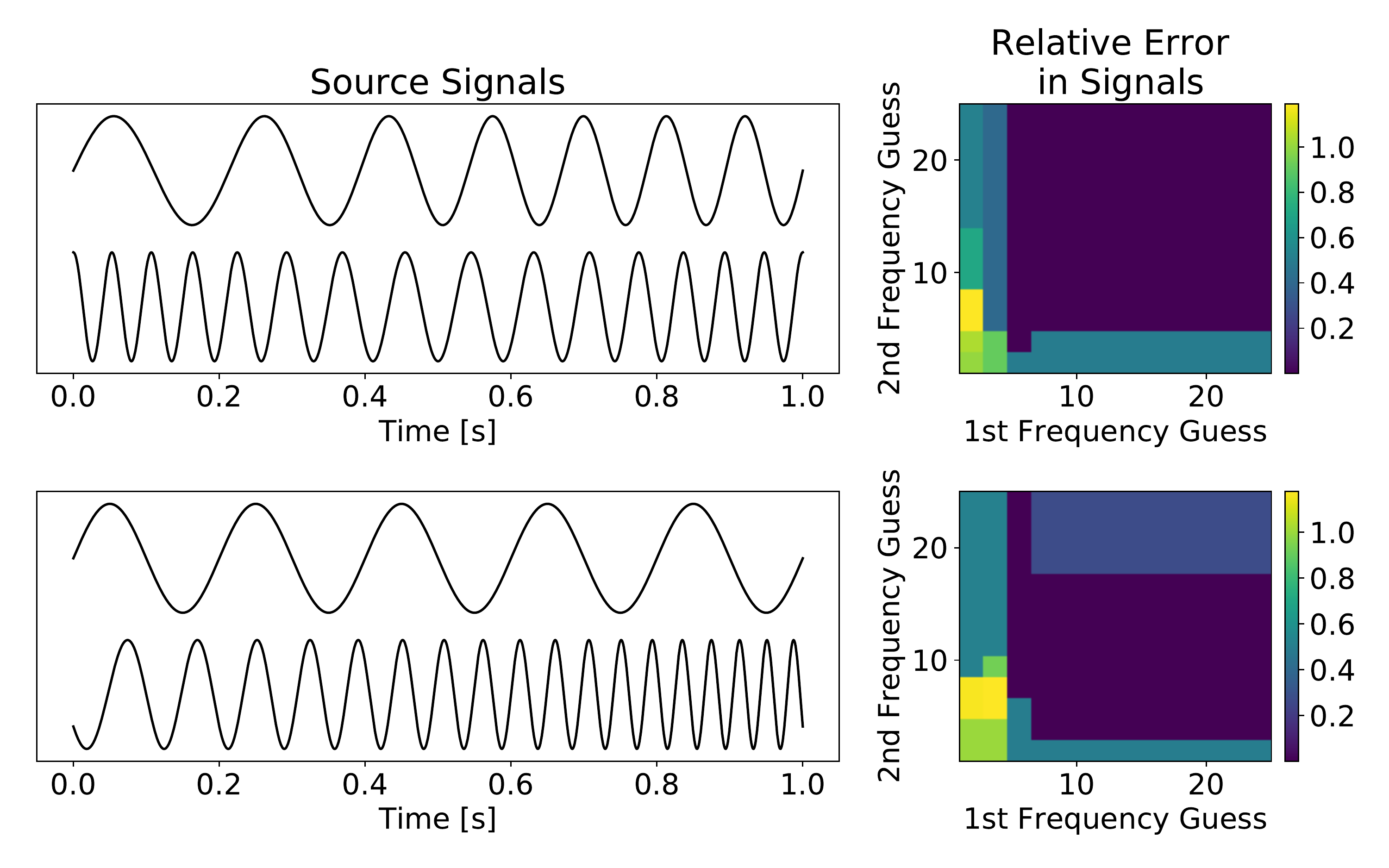}
  \caption{Two examples of the effect of initial conditions on the resulting STIMD modes. On the left are the source signals. STIMD is applied to the corresponding mixed signals. On the right is the relative error between the true source signals and the STIMD modes.}
  \label{fig:sensitivity}
\end{figure}

In addition, choosing the order in which the guesses are applied can affect the IMF's extracted from STIMD. 
Here, two examples are presented in Figure~\ref{fig:sensitivity}. 
For this system, the source signals (shown in black) are
\begin{align*}
  s_1(t) =& \sin(14 \pi t-5 \sin{(\pi t)}) \\
  s_2(t) =& \cos(30 \pi t+4 \sin(2 \pi t)),
\end{align*} 
and a $2 \times 2$ mixing matrix is used. 
STIMD is then applied to the mixed signals initial guesses for frequencies ranging from 1 $s^{-1}$ to 25 $s^{-1}$ (or equivalently $\theta'$ ranging from $2 \pi$ rad/s to $50 \pi$ rad/s).

The squared error between the STIMD modes and the true source signals is visualized on the right. 
The purple region corresponds to the region where the modes are visually correct. 
The blue regions with relative errors near $0.5$ typically correspond to when one of the two modes was correct. 
When the relative error is near $1$ (corresponding to the green and yellow regions) the modes extracted are completely incorrect.

For the second example, the source signals (like in Section \ref{subsubsec:2d}) are
\begin{align*}
s_1(t) =& \sin(10 \pi t) \\
s_2(t) =& \sin(20 \pi (t+ 0.4)^2).
\end{align*} 

It's important to note that, in both examples, the guesses need to be within a few Hz of the central frequencies of the source signals to obtain accurate results. 
Also, note that the distributions are not symmetric. 
In other words, the order in which the guesses are made for each frequency matters. 
In general, we recommend guessing frequencies in ascending order. 
This makes sense since the NMP algorithm is based on successive applications of low-pass filters.

\begin{figure}[t]
  \centering
  \begin{overpic}[scale=0.35,unit=1bp]{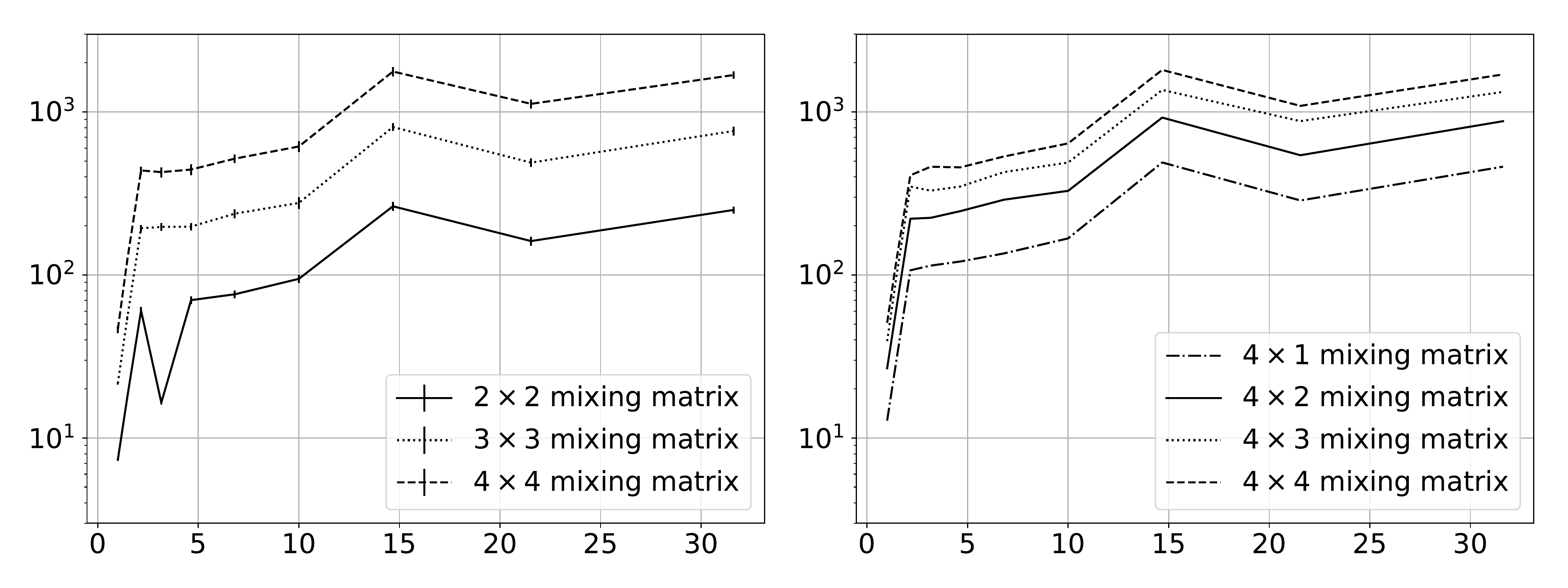}
  \put(6,-1){{Number of samples in time series}}
  \put(20,-4){[thousands]}
  \put(-1,11){\rotatebox{90}{Runtime[s]}}

  \put(56,-1){{Number of samples in time series}}
  \put(70,-4){[thousands]}
  \end{overpic}
  \vspace{0.7cm}
  \caption{Runtimes for measured for different numbers of time samples and for different mixing matrix sizes.}
  \label{fig:runtimes}
\end{figure}
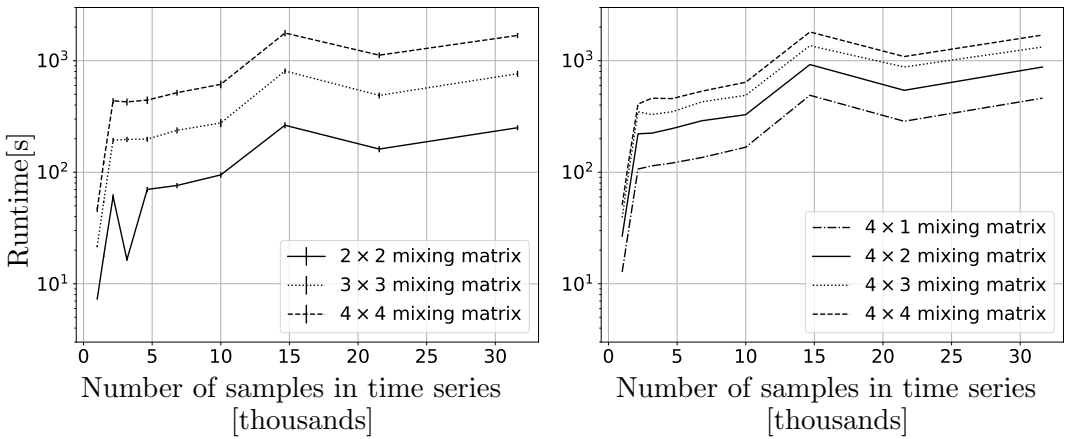

To give some indication of the computational expense of executing STIMD in practice, Figure~\ref{fig:runtimes} evaluates the run times for different number of time samples and mixing matrices. Simulations and timing scores were produced by a 32 core Intel Xeon E5-2620v4 computer with 128 GB RAM.

\section{Experiments on Real World Data}
\label{sec:demos}

In this section, we present results of STIMD on two real-world datasets in diverse domains.

\subsection{Gravitational Waves from the LIGO Experiment}
The Laser Interferometer Gravitational Wave Observatory (LIGO) is a recent Nobel prize-winning physics experiment with the goal of discovering and studying gravitational waves resulting from merging black holes~\cite{abbott2016observation}. 
The experiment consists of two detectors, one in Hanford, Washington and one in Livingston, Louisiana. 
These two detectors perform independent measurements, which can then be combined to increase confidence that a gravitational wave has been detected. 
These waves tend to be sinusoidal in nature, containing both frequency and amplitude modulation. 
The frequency modulation comes from the fact that as the black holes merge they rotate around each other with increasing frequency.

The signals measured in the first gravitational wave detection are shown in Figure \ref{fig:sensitivity}a. 
For their analysis, the LIGO collaboration computed a spectrogram (reproduced in Figure \ref{fig:sensitivity}c). 
Here the chirp corresponding to the signal is readily apparent. 
In addition, there are clearly many other residual effects from using the Fourier transform. 
For example, there are clearly nonphysical high frequency components during earlier times. 

Using initial guesses for the phases of $\theta_1(t)/2 \pi = 50 t$ and $\theta_2(t)/2 \pi = 128 t$, we obtained the IMF shown in Figure \ref{fig:sensitivity}b, which clearly corresponds to the primary chirp seen in the data. 

\begin{figure}[t]
\centering
\begin{overpic}[scale=0.35,unit=1bp]{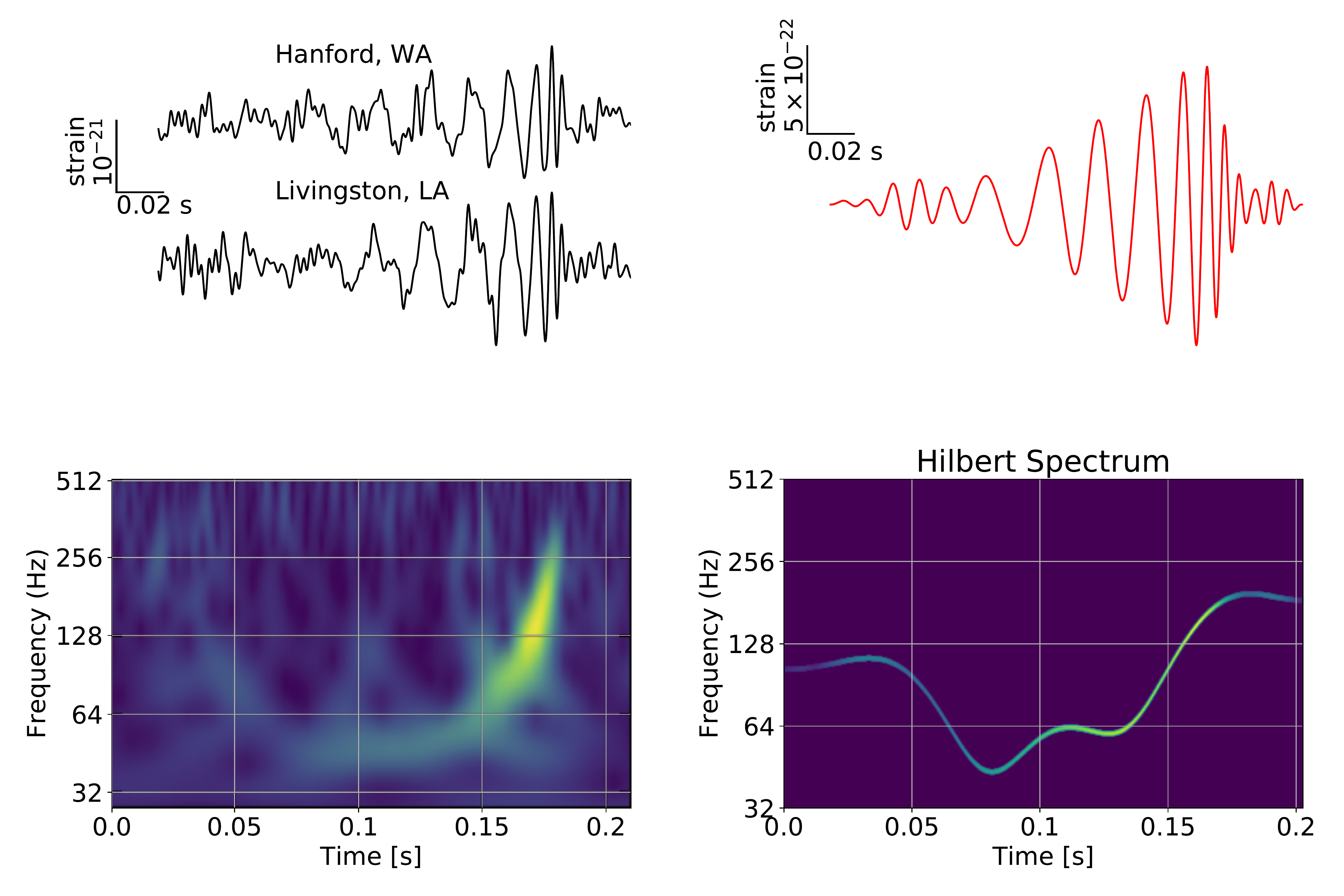}
\put(10,62){(a)}
\put(50,62){(b)}
\put(10,35){(c)}
\put(60,35){(d)}
\end{overpic}
\caption{Application of STIMD to sample data from LIGO experiment. (a) Signals measured in LIGO experiment using two detectors. (b) STIMD modes. (c) Spectrogram computed in LIGO analysis. (d) Hilbert spectrum computed using STIMD modes. Yellow corresponds to frequencies of greater intensity while purple corresponds to frequencies of lower intensity.}
\label{fig:ligo}
\end{figure}

\subsection{Neural Recordings from Rodent Hippocampus}

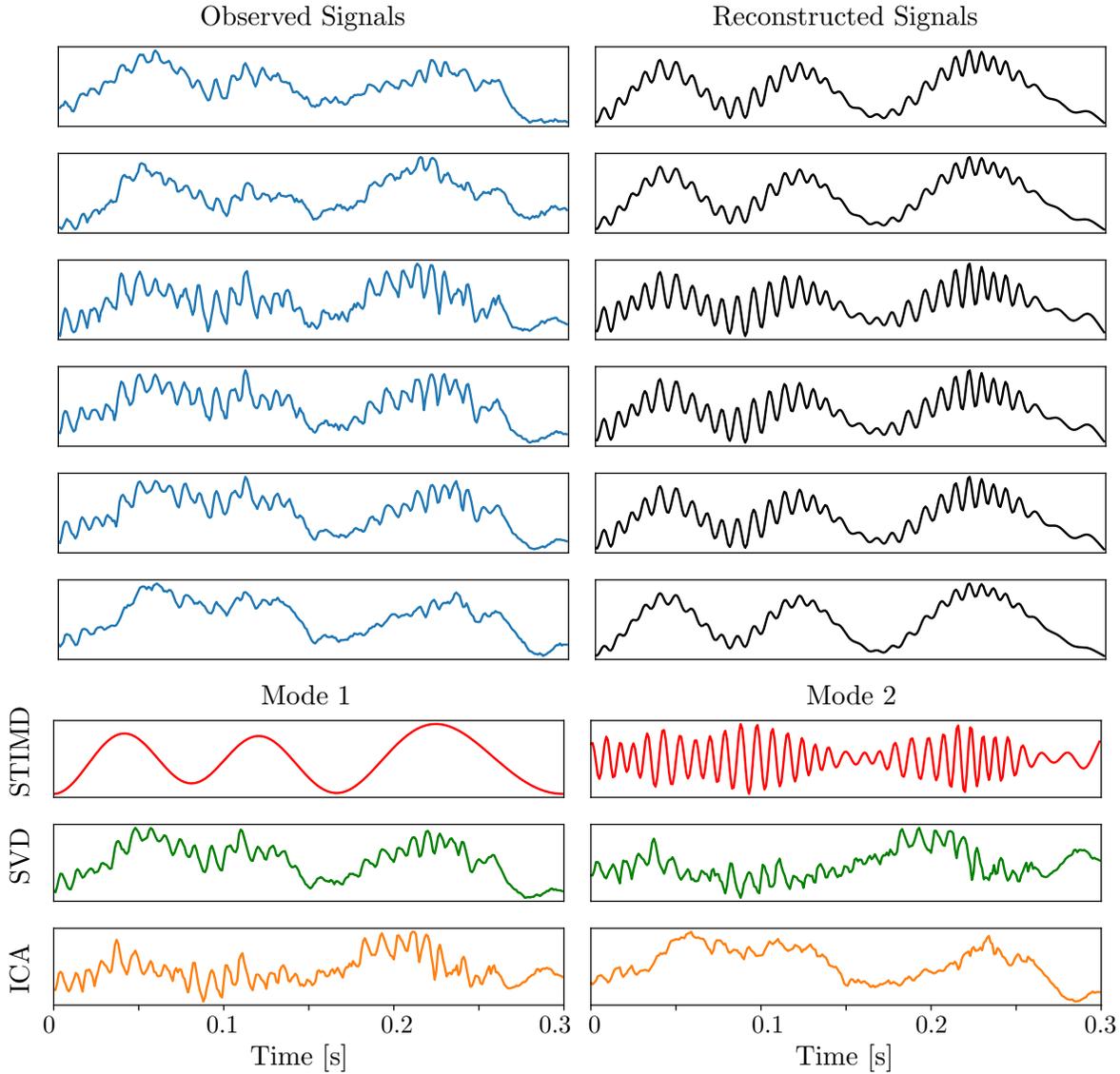
\begin{figure}[t]
\centering
\begin{overpic}[scale=0.6,unit=1bp]{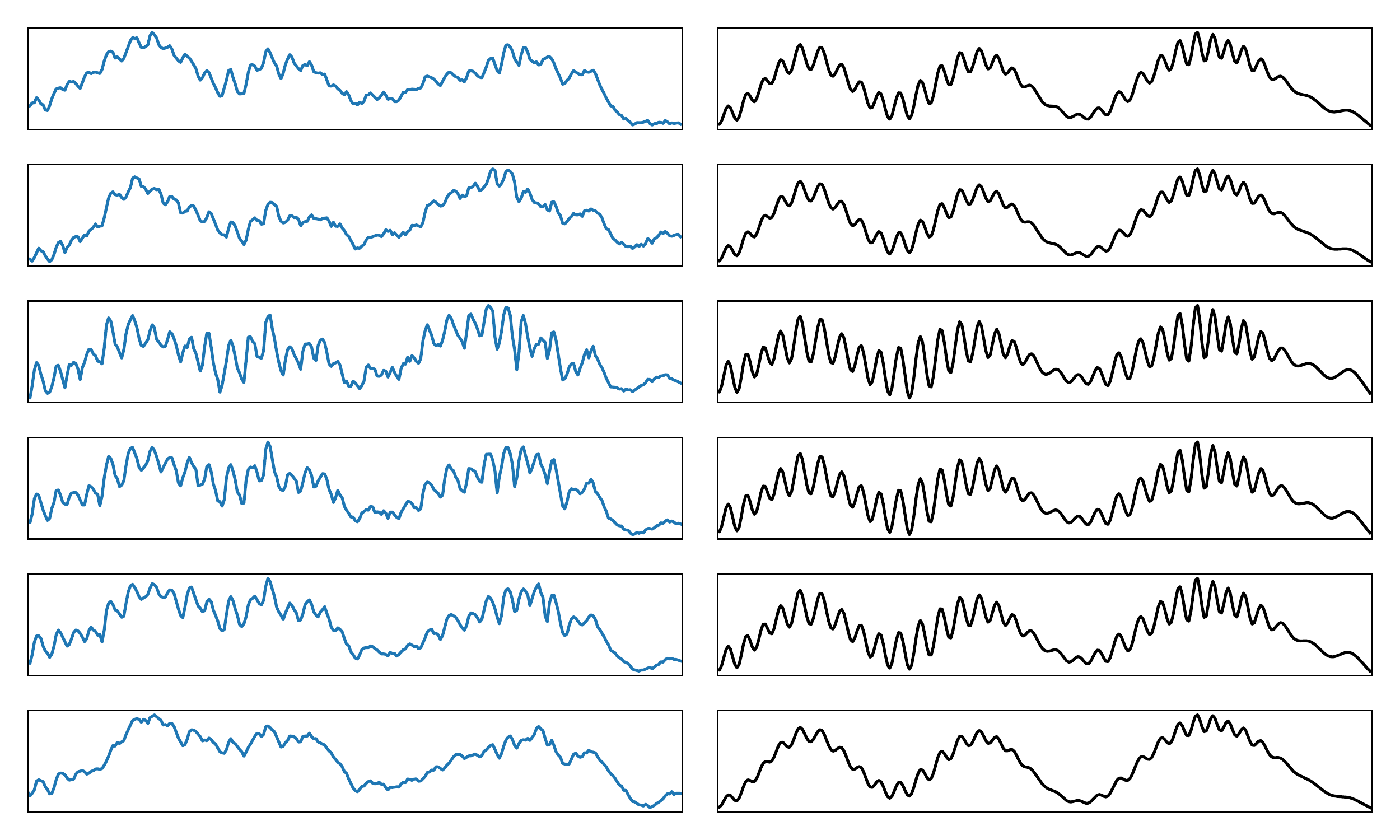}
\put(15,60){Observed Signals}
\put(62,60){Reconstructed Signals}
\end{overpic}
\par\medskip 
\begin{overpic}[scale=0.6,unit=1bp]{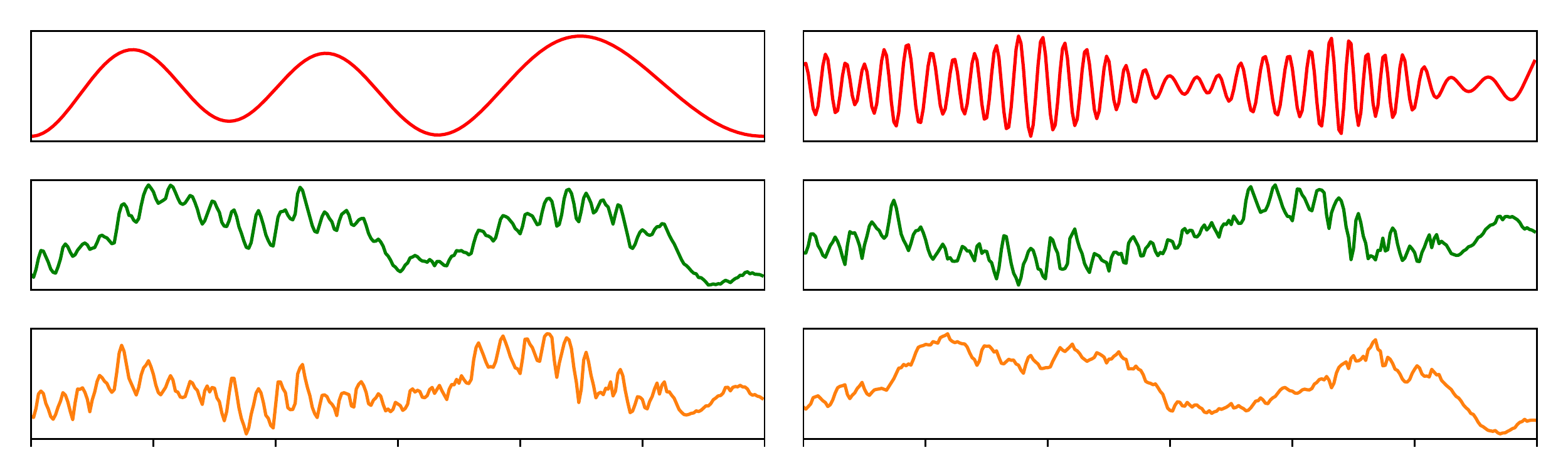}
\put(21,29.5){Mode 1}
\put(71,29.5){Mode 2}
\put(-2,21){\rotatebox{90}{STIMD}}
\put(-2,13){\rotatebox{90}{SVD}}
\put(-2,3){\rotatebox{90}{ICA}}

\put(20,-3.5){Time [s]}
\put(70,-3.5){Time [s]}

\put(1,-0.5){\small 0}
\put(16,-0.5){\small 0.1}
\put(32,-0.5){\small 0.2}
\put(46.5,-0.5){\small 0.3}

\put(51,-0.5){\small 0}
\put(66,-0.5){\small 0.1}
\put(81,-0.5){\small 0.2}
\put(96.5,-0.5){\small 0.3}

\end{overpic}
\vspace*{7mm}
\caption{Application of STIMD to sample rodent neural recordings in hippocampus. Blue: Measured neural recordings. Red: STIMD modes. Green: SVD modes. Orange: ICA modes. Black: Signals reconstructed using STIMD modes and mixing matrix.} 
  \label{fig:hippocampus}
\end{figure}

As a second example, we analyze recordings of neural activity. 
This data is available at \cite{hippocampus-data} and has been described previously~\cite{grosmark2016diversity}. 
It is well known that these local field potential signals in the rodent hippocampus contain rhythmic activity, with frequency and amplitude modulation~\cite{nguyen2008instantaneous}. 

Here we analyze a subset of the recordings from six neighboring electrodes, placed at equal spacing, over a 300ms period. 
Using peaks in the signals, we chose initial guesses for the phases $\theta_1(t) = 20 \pi t - \pi / 2$ and $\theta_2(t) = 200 \pi t$, respectively. 
The mixing matrix $B$ consequently has dimensions $6 \times 2$. 
In Figure \ref{fig:hippocampus}, we show the observed signals (blue), the STIMD modes (red), the SVD modes (green), the ICA modes (orange), and the reconstructed signals (black) computed by multiplying the mixing matrix B and by the STIMD modes. 
Note that in contrast to ICA and SVD, STIMD is able to separate the modes by their frequencies. 
It is unclear if the SVD or ICA decompositions are interpretable for this case. 
In addition, STIMD extracts spatial modes; in particular, the second STIMD mode, which corresponds to a high frequency wave, has higher amplitude in the central electrodes.


\section{Conclusions}
\label{sec:conclusions}

Principled mathematical methods for spatiotemporal decompositions are critically enabling for many emerging large-scale applications across the physical and biological sciences.  
Of particular interest is the ability to perform blind source separation on data generated from nonlinear and non-stationary dynamical processes.  
Our proposed STIMD mathematical architecture provides a compromise between the commonly used methods of ICA, SVD, and DMD.  
Specifically, STIMD is not constrained like DMD to model stationary signals and Fourier modes in time.  
However, we still make use of the constraint that our data must obey a certain set of dynamics commonly found in physical systems. 
In addition, we have the ability to compute a Hilbert spectrum and perform future state prediction, which is not guaranteed for ICA or SVD.   

In this paper, the STIMD method is applied on synthetic data to evaluate its feature extraction performance.  The method leverages recent key innovations for signal processing from EMD and NMP.  
Indeed, by exploiting the IMF time constraint, the STIMD method frames an optimization problem that extracts meaningful features and low-rank modes from spatiotemporal data.  
We demonstrate the method on two real-world data sets, the LIGO experiment for the discovery of gravity waves and neural activity recordings from the rodent hippocampus.  
In both cases, the STIMD method produces a clean spatiotemporal decomposition with interpretable modes.  
We suggest STIMD as a method for data-driven discovery that may be widely applied to many domains with spatiotemporal data.

\appendix
\section{NMP Minimization Implementation} \label{app:nmp}
For the NMP method we must solve the minimization problem
\begin{equation}
  \min_{a} \abs{r_k - a(\theta(t)) \cos{(\theta}(t))}_2^2 \text{ where } a(\theta) \in V(\theta).
\end{equation} 
In Ref.~\cite{hou2013data}, Shi and Hou solve this by noting that constraining $a(\theta)$ to lie in $V(\theta)$ is equivalent to applying a low-pass filter in $\theta$-space.  The corresponding algorithm is shown in 

\begin{algorithm}
\caption{Minimization Algorithm for Nonlinear Matching Pursuit (NMP) Method for Periodic Data}
\begin{algorithmic}[1]
\STATE Input: measured signal $x(t)$ and phase function of IMF's $\theta(t)$.  
\STATE Output: $a(t),b(t)$
\STATE Define the normalized phase function $\bar{\theta}(t) = \frac{\theta(t) - \theta(0)}{\theta(T) - \theta(0)}$
\STATE{$L_{\theta} = \frac{\theta(T) - \theta(0)}{2 \pi}$}

\STATE{$x(\theta) :=  \text{Interpolate}(x(t),\theta(t))$ \COMMENT{Reexpress $x(t)$ in terms of the $\theta$ coordinate}}.
\STATE{$\hat{x}(\omega) = \mathcal{F}\left( x(\theta) \right)$} 
\STATE{$a(\theta) := \mathcal{F}^{-1} \left( \hat{x}(\omega + L_{\theta}) + \hat{x}(\omega - L_{\theta}) \cdot \chi_{\lambda}(\omega / L_{\theta}) \right)$}
\STATE{$b(\theta) := \mathcal{F}^{-1} \left( i \cdot \left( \hat{x}(\omega + L_{\theta}) - \hat{x}(\omega - L_{\theta}) \cdot \chi_{\lambda}(\omega / L_{\theta}) \right) \right) $}
\STATE{$a(t) = \text{Interpolate}(a(\theta),t)$}
\STATE{$b(t) = \text{Interpolate}(b(\theta),t)$}
\RETURN $a,b$
\end{algorithmic}
\label{alg:NMP2}
\end{algorithm}
In this algorithm $\mathcal{F}$, and $\mathcal{F}^{-1}$ denote the Fourier Transform, and inverse Fourier Transform respectively 
\begin{gather*}
  \mathcal{F}(r) = \frac{1}{N} \sum_{j = 1}^M r_j e^{-i 2 \pi \omega \bar{\theta}_j}, \omega = -N/2 + 1, \cdots, N/2 \\
  \mathcal{F}^{-1}(\hat{r}) = \frac{1}{N} \sum_{\omega = -N/2 + 1}^{N/2} \hat{r} e^{i 2 \pi \omega \bar{\theta}_{k,j}^n}, j = 0,\ldots, N-1
\end{gather*}
$\chi_{\lambda}(\omega)$ is the cutoff function used in the low-pass filter. Here we use the function 
\[\chi_{\lambda}(\omega) = \begin{cases} 
      1 + \cos(\pi \omega/\lambda) & -\lambda < \omega < \lambda \\
      0 & \text{otherwise}
   \end{cases}
\]

\section*{Acknowledgments}
We are grateful for discussions with K. D. Harris and T. S. Tan. In addition, we thank T. Y. Hou and Z. Shi for their support and providing us with the corresponding MATLAB versions of their NMP code.
This work was funded by the AFOSR grant (FA9550-17-1-0329) to JNK; the NSF (award 1514556), the Alfred P. Sloan Foundation, and the Washington Research Foundation to BWB.

\end{document}